\documentclass{autart}

\RequirePackage{fix-cm}

\usepackage{graphicx}
\usepackage{psfrag}
\usepackage{caption}

\usepackage{subfig}

\usepackage{url}
\usepackage{color}
\usepackage{cite}
\usepackage{epsfig}
\usepackage{hyperref}
\usepackage{amsmath,amssymb}
\usepackage{slashbox}
\usepackage{tablefootnote}
\usepackage{longtable}

\newcommand{\rset}{\mathbb{R}}

\newcommand{\Ncal}{\mathcal{N}}

\newcommand{\lu}{\mathcal{L}}
\newtheorem{theorem}{Theorem}[section]
\newtheorem{assumption}[theorem]{Assumption}

\newtheorem{lemma}[theorem]{Lemma}

\usecounter{algorithmenumi}

\begin{document}

\begin{frontmatter}

\date{October 1, 2013}

\title{On linear convergence of a distributed dual gradient algorithm
for linearly constrained separable  convex problems}

\author{Ion Necoara and Valentin Nedelcu}
\address{Automatic Control and Systems Engineering Department, University Politehnica Bucharest, Romania.}

\thanks{Corresponding author: Ion Necoara, email: ion.necoara@acse.pub.ro.}

\begin{abstract}
In this paper we propose a distributed dual gradient  algorithm for
minimizing  linearly constrained separable  convex problems and
analyze its rate of convergence. In particular, we prove that under
the assumption of  strong convexity and Lipshitz continuity of the
gradient of the primal objective function we have a global error
bound type property for the dual problem. Using this error bound
property we devise a  fully distributed dual gradient scheme, i.e. a
gradient  scheme based on a weighted step size,  for which we derive
global linear rate of convergence for both dual and primal
suboptimality and for primal feasibility violation.  Many real
applications, e.g. distributed model predictive control, network
utility maximization or optimal power flow, can be posed as linearly
constrained separable  convex problems for which dual gradient type
methods from literature  have sublinear convergence rate. In the
present paper we prove for the first time that in fact we can
achieve linear convergence rate for such algorithms when they are
used for solving these applications.  Numerical simulations are also
provided to confirm our theory.
\end{abstract}

\end{frontmatter}


\section{Introduction}
\label{sec_introduction} Nowadays, many engineering applications
which appear in the context of communications networks or networked
systems can be posed as large scale linearly constrained  separable convex
problems.  Several important applications that can be modeled in
this framework, distributed model predictive control (DMPC) problem
for networked systems \cite{NecNed:11},  the network utility
maximization (NUM) problem \cite{BecNed:13}, and the direct current optimal power
flow (DC-OPF) problem for a power system \cite{BakBis:03},   have
attracted great attention lately. Due to the large dimension and the
separable structure of these problems, distributed optimization
methods have become an appropriate tool for solving them.

Distributed optimization methods are based on decomposition
\cite{NecNed:11}. Decomposition methods represent a powerful tool
for solving these types of problems due to their ability of dividing
the original large scale problem into smaller subproblems which are
coordinated by a master problem. Decomposition methods can be
divided into two main classes: primal and dual decomposition. While
in the primal decomposition methods the optimization problem is
solved using the original formulation and variables
\cite{FarSca:12,NecNed:11},   in dual decomposition the constraints
are moved into the cost using the Lagrange multipliers and then the
dual problem is solved \cite{SenShe:86,NecSuy:08}. In many
applications, such as (DMPC), (NUM) and (DC-OPF) problems, when the
constraints set is complicated (i.e. the projection on this set is
hard to compute), dual decomposition is more effective since a
primal approach would require at each iteration a projection onto
the feasible set, operation that is numerically  expensive.

First order decomposition methods for solving dual problems have
been extensively studied in the literature. Dual subgradient methods
based on averaging, that produce primal solutions in the limit, can
be found e.g. in \cite{SenShe:86}. Convergence rate analysis for the
dual subgradient method is given e.g. in \cite{NedOzd:09}, where
estimates of order $\mathcal{O} (1/\sqrt{k})$ for suboptimality and
feasibility violation of an average primal sequence are provided,
with $k$ denoting the iteration counter. In \cite{NecSuy:08} the
authors derive a dual decomposition method based on  a  fast
gradient algorithm and  a smoothing  technique  and prove rate of
convergence of order $\mathcal{O}\left(\frac{1}{k}\right)$ for
primal suboptimality and feasibility violation for an average primal
sequence. Also, in \cite{NecNed:13,NedNec:12} (\cite{PatBem:12}) the
authors propose inexact (exact) dual fast gradient algorithms for
which estimates of order $\mathcal{O}\left(\frac{1}{k^2}\right)$ in
an average primal sequence  are provided for  primal  suboptimality
and feasibility violation. For the special case of QP problems, dual
fast gradient algorithms were also analyzed in
\cite{GisDoa:13,KogFin:11}. To our knowledge, the first result on
the  linear convergence of dual gradient methods was provided in
\cite{LuoTse:93a}. However, the authors in \cite{LuoTse:93a} were
able to show linear convergence only \textit{locally} using a local
error bound condition that estimates the distance from the dual
optimal solution set in terms of norm of a proximal residual.
Another strand of this literature uses alternating direction method
of multipliers (ADMM) \cite{HonLuo:12,WeiOzd:13} or  Newton methods
\cite{WeiOzd:10,NecSuy:09}. For example, \cite{HonLuo:12}
established a  linear convergence rate of (ADMM) using an error
bound condition that holds under specific assumptions on the primal
problem, while in \cite{WeiOzd:13} sublinear rate of convergence is
proved for (ADMM), but for more general assumptions on the primal
objective function. In \cite{WeiOzd:10,NecSuy:09} distributed Newton
algorithms are derived with fast convergence under the assumption
that the primal objective function is self-concordant. Finally, very
few results were known in the literature on distributed
implementations of dual gradient type methods since most of the
papers enumerated above require a centralized  step size. Recently,
in \cite{NecCli:13} (\!\cite{BecNed:13}), distributed (dual fast)
gradient algorithms are given,  where the step size is chosen
distributively, and estimates of order
$\mathcal{O}\left(\frac{1}{k}\right)$  for primal suboptimality and
infeasibility  in the last primal iterate (linear) are given.

Despite widespread use of gradient methods for solving dual
problems, there are some aspects that have not been fully studied.
In particular, in applications the main interest is in finding
primal vectors that are near-feasible and near-optimal. We also need
to devise algorithms with fast convergence rate, e.g. linear
convergence. Finally, we are interested in providing distributed
schemes, i.e. methods based on distributed computations. These
represent the main issues that we pursue in this paper.

\textit{Contributions:} In this paper we propose a  distributed dual
gradient method generating approximate primal feasible and optimal
solutions but with great improvement on the convergence rate w.r.t.
existing results. Under the assumptions of strong convexity and
Lipschitz continuity of the gradient of the primal objective
function, which are often satisfied in practical applications (e.g.
(DMPC), (NUM) or (DC-OPF)), we prove that the corresponding dual
problem satisfies  a certain global error bound property that
estimates the distance from the dual optimal solution set in terms
of the norm of a proximal residual. In order to prove such a result
we extend the approach developed in \cite{WanLin:13,LuoTse:93a},
where the authors show such a property for objective functions
having a certain structure and  constraints set defined in  terms of
bounded polyhedra, to the case where the objective function is more
general and  the constraints set is an unbounded polyhedron.  This
nontrivial extension also allows us to tackle dual problems, where
e.g. the constraints are defined in terms of the nonnegative
orthant. In these settings we analyze the convergence behavior of a
distributed dual gradient algorithm, for which we are able to
provide for the first time \textit{global} linear convergence rate
on primal suboptimality and feasibility violation for the last
primal iterate, as opposed  to the results in \cite{LuoTse:93a}
where only \textit{local} linear convergence was derived for such an
algorithm. Moreover, our algorithm is fully distributed since is
based on a weighted step size, as opposed to typical dual
distributed schemes existing in literature, where a centralized step
size is used and sublinear convergence is proved
\cite{DoaKev:11,GisDoa:13,NecNed:13,MeiUlb:13}. Note that our
results are also related to those in \cite{HonLuo:12}: in
particular, paper \cite{HonLuo:12} established  an error bound
property for the \textit{augmented} dual function and then proved
linear convergence for the (ADMM) method. However, the main
drawbacks with (ADMM) consist of  the difficulty in tuning the
penalty parameter in the augmented Lagrangian and the centralized
choice~of~it.

\noindent \textit{Paper Outline:} In Section \ref{sec_formulation}
we introduce our optimization model and discuss the (DMPC) problem
for a networked system. In Section \ref{sec_error_bound} we prove a
certain error bound property of the dual function which allows  us
to derive global linear converge for a fully distributed dual
gradient method in Section \ref{sec_dg_error_bound}. In Sections
\ref{sec_distributed}  we discuss implementation issues, in
particular  in the context of (DMPC), and finally in Section
\ref{sec_numerical} we provide some numerical simulations that
confirm our theory.

\noindent \textit{Notations:}  For $z, y \in \rset^n$ we denote the
Euclidean inner product $\langle z,y \rangle= z^Ty= \sum_{i=1}^n z_i
y_i$, the Euclidean norm  $\left \| z \right \|=\sqrt{\langle z, z
\rangle}$ and the infinity norm  $\|z\|_{\infty}=\sup_i |z_i|$. For
a matrix $G \in \rset^{m \times n}$, $\|G\|$ denotes its spectral
norm. Also, we denote the orthogonal  projection onto the
nonnegative orthant $\rset^n_+$ by $\left[ z \right]_+$ and the
orthogonal projection onto the convex set $D$ by  $[z]_D$. For a
positive definite matrix $W$ we define norm $\|z\|_W=\sqrt{z^TWz}$
and the projection of  vector $z$ onto a convex set $D$ w.r.t. norm
$\|\cdot\|_W$ by $[z]_D^W$. For a matrix $A$,  $A_{i}$ is its $i$th
(block)~column.


\section{Problem formulation}
\label{sec_formulation}

We consider the following large scale linearly constrained separable convex
optimization problem:
\begin{align}
\label{eq_prob_princc} f^* = &\min_{z_i \in \rset^{n_i}}
f(z)~~~\left(= \sum_{i=1}^M f_i(z_i)\right)\\
&\text{s.t.:} ~~Az=b, ~~Cz\leq c, \nonumber
\end{align}
where $f_i: \rset^{n_i} \to \rset^{}$ are  convex functions,
$z=\left[z_1^T \cdots z_M^T\right]^T \in \rset^n$, $A \in \rset^{p
\times n}$, $C \in \rset^{q \times n}$, $b\in \rset^p$ and $c \in
\rset^q$. To our problem \eqref{eq_prob_princc} we associate a
bipartite communication graph $\mathcal{G}=\left(V_1,V_2,E\right)$,
where $V_1=\left\{1,\dots,M\right\}$,
$V_2=\left\{1,\dots,\bar{M}\right\}$ and $E \in
\left\{0,1\right\}^{\bar{M} \times M}$ represents the incidence
matrix. E.g., in the context of (NUM) and (DC-OPF), $V_1$ denotes
the set of sources, $V_2$ the set of links between sources and the
incidence matrix $E$ models the way sources interact. In (DMPC),
$V_1=V_2$ represents the set of interacting subsystems, while the
incidence matrix $E$ indicates the dynamic couplings between these
subsystems. We assume that $A$ and $C$ are block matrices with the
blocks $A_{ji} \in \rset^{p_j \times n_i}$ and $C_{ji} \in
\rset^{q_j \times n_i}$, where $\sum_{i=1}^M n_i=n$,
$\sum_{j=1}^{\bar{M}} p_j=p$ and $\sum_{j=1}^{\bar{M}} q_j=q$. We
also assume that if $E_{ji}=0$, then both blocks $A_{ji}$ and
$C_{ji}$ are zero. In these settings we allow a block $A_{ji}$ or
$C_{ji}$ to be zero even if $E_{ji}=1$. We also introduce the index
sets:
\[ \bar{\Ncal}_i=\left\{j \in V_2 : E_{ji} \neq 0\right\} \;\;   \text{and} \;\;
\Ncal_j = \left\{i \in V_1 : E_{ji} \neq 0\right\} \] for all $i
\in V_1, j \in V_2$, which describe the local information flow in the graph. Note
that the cardinality of the sets $\bar{\Ncal}_i$ and $\Ncal_j$ can
be viewed as a measure for the degree of separability of problem
\eqref{eq_prob_princc}. Therefore, the local information structure
imposed by the graph $\mathcal{G}$ should be considered as part of
the problem formulation. Further, we make the following assumptions
regarding the optimization problem \eqref{eq_prob_princc}:
\begin{assumption}
\label{ass_strong} $(a)~$ The functions $f_i$ have Lipschitz
continuous gradient with constants $L_i$ and are $\sigma_i$-strongly
convex w.r.t. the Euclidean
norm $\|\cdot\|$ on $\rset^{n_i}$ \cite{Nes:04}.\\
$(b)~$ Matrix $A$ has full row rank and there
exists a feasible point $\tilde z$ for problem
\eqref{eq_prob_princc}  such that $A \tilde z = b$ and  $C \tilde z < c$.
\end{assumption}

Note that if Assumption \ref{ass_strong} $(a)$ does not hold, we can
apply smoothing techniques by adding a regularization term to the
function $f_i$ in order to obtain a strongly convex approximation of
it (see e.g. \cite{NecSuy:08} for more details).
 Assumption \ref{ass_strong} $(b)$ implies that
strong duality holds for optimization problem \eqref{eq_prob_princc}
and the set of optimal Lagrange multipliers is bounded \cite[Theorem
2.3.2]{HirLem:96}. Note that Assumption \ref{ass_strong} $(b)$ is
not restrictive: we can always remove the redundant equalities so
that matrix $A$ has full row rank and strict feasibility for the
inequality constraints is usually satisfied in applications (e.g.
(DMPC), (NUM) or (DC-OPF)). In particular, we have:
\begin{equation}
\label{eq_dual_prob} f^*=\max_{\nu \in \rset^p,\mu \in \rset^q_+}
d(\nu,\mu),
\end{equation}
where $d(\nu,\mu)$ denotes the dual function of
\eqref{eq_prob_princc}:
\begin{equation}
\label{eq_dual_func} d(\nu,\mu)=\min_{z \in \rset^n} \lu
(z,\nu,\mu),
\end{equation}
with  the Lagrangian function \[\lu(z,\nu,\mu) = f(z)+\langle \nu,
Az-b \rangle +\langle \mu, Cz-c \rangle.\]
For simplicity of the exposition we introduce further the following notations:
\begin{equation}
\label{eq_notations} G=\left[
\begin{array}{c}
  A \\
  C \\
\end{array}
\right]~\text{and} ~ g=\left[
\begin{array}{c}
  b \\
  c \\
\end{array}
\right].
\end{equation}
Since $f_i$ are strongly convex functions, then $f$ is also strongly
convex w.r.t. the Euclidean norm $\| \cdot\|$ on $\rset^n$, with
convexity parameter e.g.  $\sigma_{\mathrm{f}}=\min
\limits_{i=1,\dots,M} \sigma_i$. Further,  the dual function $d$ is
differentiable and its gradient is given by the following expression
\cite{NecNed:13}:
\begin{equation*}
\nabla d(\nu,\mu)=Gz(\nu,\mu)-g,
\end{equation*}
where $z(\nu,\mu)$ denotes the unique optimal solution of the
inner problem  \eqref{eq_dual_func}, i.e.:
\begin{equation}
\label{eq_inner_sol} z(\nu,\mu)=\arg\min_{z \in \rset^n}
\lu(z,\nu,\mu).
\end{equation}
Moreover, the gradient $\nabla d$ of the dual function is Lipschitz
continuous w.r.t. Euclidean  norm $\|\cdot\|$, with constant
\cite{NecNed:13}:
\begin{equation*}
L_{\text{d}}=\frac{\left\|G\right\|^2}{\sigma_{\mathrm{f}}}.
\end{equation*}
If we denote by $\nu_{\bar{\Ncal}_i}=\left[\nu_j\right]_{j \in \bar{\Ncal}_i}$ and
by $\mu_{\bar{\Ncal}_i}=\left[\mu_j\right]_{j \in \bar{\Ncal}_i}$ we can observe
that the dual function can be written in the following separable
form:
\begin{equation*}
d(\nu,\mu)=\sum_{i=1}^M d_i(\nu_{\bar{\Ncal}_i},\mu_{\bar{\Ncal}_i})-\langle
\nu,b\rangle-\langle \mu, c\rangle,
\end{equation*}
with
\begin{align}
\label{eq_separable_dual}
d_i(\nu_{\bar{\Ncal}_i},\mu_{\bar{\Ncal}_i})&=\min_{z_i \in \rset^{n_i}}
f_i(z_i)+ \langle
\nu, A_{i} z_i\rangle+\langle \mu, C_{i} z_i\rangle\\
&=\min_{z_i \in \rset^{n_i}} f_i(z_i)+ \sum_{j \in \bar{\Ncal}_i}
\left\langle A_{ji}^T \nu_j+C_{ji}^T\mu_j,z_i\right\rangle.
\nonumber
\end{align}
In these settings, we have that the gradient $\nabla d_i$ is:
\begin{equation*}
\nabla d_i(\nu_{\bar{\Ncal}_i},\mu_{\bar{\Ncal}_i})=\left[\begin{array}{c}
  \left[A_{ji}\right]_{j \in
\bar{\Ncal}_i} \\
  \left[C_{ji}\right]_{j \in
\bar{\Ncal}_i} \\
\end{array}\right] z_i(\nu_{\bar{\Ncal}_i},\mu_{\bar{\Ncal}_i}),
\end{equation*}
where $z_i(\nu_{\bar{\Ncal}_i},\mu_{\bar{\Ncal}_i})$ denotes the
unique optimal solution in \eqref{eq_separable_dual}. Note that
$\nabla d_i$ is Lipschitz continuous  w.r.t. Euclidean  norm $\|
\cdot \|$,  with constant \cite{NecNed:13}:
\begin{equation}
L_{\text{d}_i}=\frac{\left\|\left[
\begin{array}{c}
  \left[A_{ji}\right]_{j \in
\bar{\Ncal}_i} \\
   \left[C_{ji}\right]_{j \in
\bar{\Ncal}_i} \\
\end{array}
\right]\right\|^2}{\sigma_i}.
\end{equation}
For simplicity of the exposition we will consider the
notation $\lambda=\left[\nu^T \mu^T \right]^T$ and we will also denote the effective  domain of the dual function by $\mathbb{D}=\rset^p \times \rset^q_+$. The following result,
which is a distributed version of descent lemma is central in our
derivations of a distributed dual algorithm and in the proofs of its convergence rate.

\begin{lemma}
\label{descent_lemma} Let Assumption \ref{ass_strong} $(a)$ hold.
Then, the following inequality is valid:
\begin{equation}
\label{eq_descent} d(\lambda) \!\geq\! d(\bar{\lambda})+\left\langle
\nabla d(\bar{\lambda}), \lambda\!-\!\bar{\lambda}\right\rangle
-\frac{1}{2}\|\lambda-\!\bar{\lambda}\|_W^2 ~~ \forall \lambda,
\bar{\lambda} \in \mathbb{D},
\end{equation}
where the matrix $W=\text{diag}(W_{\nu},W_{\mu})$ with the matrices
$W_{\nu}=\text{diag}\left(\sum_{i \in \Ncal_j
}L_{\text{d}_i}I_{p_j}; j \in V_2\right)$ and
$W_{\mu}=\text{diag}\left(\sum_{i \in \Ncal_j
}L_{\text{d}_i}I_{q_j}; j \in V_2\right)$.
\end{lemma}
\pf A similar result for the case of inequality constraints was
given in \cite{BecNed:13,NecCli:13}.  Let us first denote
$\lambda_{\bar{\Ncal}_i}=\left[\nu_{\bar{\Ncal}_i}^T~
\mu_{\bar{\Ncal}_i}^T\right]^T$. Using now the continuous Lipschitz
gradient property of $d_i$ we can write for each $i=1,\dots,M$
\cite{Nes:04}:
\begin{align*}
d_i(\lambda_{\bar{\Ncal}_i})&\geq
d_i(\bar{\lambda}_{\bar{\Ncal}_i})+\left\langle \nabla
d_i(\bar{\lambda}_{\bar{\Ncal}_i}),\lambda_{\bar{\Ncal}_i}
-\bar{\lambda}_{\bar{\Ncal}_i}\right\rangle\\
&\qquad -\frac{L_{\text{d}_i}}{2}\|\lambda_{\bar{\Ncal}_i}-\bar{\lambda}_{\bar{\Ncal}_i}\|^2~~~\forall \lambda_{\bar{\Ncal}_i},\bar{\lambda}_{\bar{\Ncal}_i}.
\end{align*}
Summing up these inequalities for all $i=1,\dots,M$ and adding
$\langle\lambda, \left[b^T c^T\right]^T\rangle$ to both sides of
the previous inequality we obtain:
\begin{align*}
d(\lambda) \geq d(\bar{\lambda}) + \left\langle \nabla
d(\bar{\lambda}), \lambda - \bar{\lambda} \right \rangle
-\sum_{i=1}^M
\frac{L_{\text{d}_i}}{2}\|\lambda_{\bar{\Ncal}_i}-\bar{\lambda}_{\bar{\Ncal}_i}\|^2.
\end{align*}
Using now the definition of $\lambda_{\bar{\Ncal}_i}$ we can
write:
\begin{align*}
& \sum_{i=1}^M
L_{\text{d}_i}\|\lambda_{\bar{\Ncal}_i}-\bar{\lambda}_{\bar{\Ncal}_i}\|^2\\
&=\sum_{i=1}^M L_{\text{d}_i} \sum_{j=1}^{\bar{M}} E_{ji}\left(\|\nu_j-\bar{\nu}_j\|^2+\|\mu_j-\bar{\mu}_j\|^2\right)\\
&=\sum_{j=1}^{\bar{M}} \left(\|\nu_j-\bar{\nu}_j\|^2+\|\mu_j-\bar{\mu}_j\|^2\right) \left(\sum_{i=1}^M L_{\text{d}_i} E_{ji}\right)
\end{align*}
Introducing this result into the previous inequality and using the definition of $W$ we conclude
the statement. \qed

\textit{Tightness of the descent lemma}. Our descent lemma (Lemma
\ref{descent_lemma}) is ``tight'' in the following sense: there are
functions for which $W$ in \eqref{eq_descent} cannot be replaced by
smaller diagonal  matrices in positive definite sense. We show this
on a simple example. Let $f_i(z_i) = \frac{\sigma_i}{2} z_i^T z_i$
and $n_i=1$. In this case $d(\lambda) = - \frac{1}{2} \lambda^T
\left( G Q^{-1} G^T \right) \lambda - g^T \lambda$, where $Q =
\text{diag}(\sigma_i; i \in V_1)$.  Note that we can write
$d(\lambda) = d(\bar{\lambda}) + \left\langle \nabla
d(\bar{\lambda}), \lambda - \bar{\lambda}\right \rangle  -
\frac{1}{2} \|\lambda- \bar{\lambda}\|_{G Q^{-1} G^T }^2$. Let us
define $\lambda- \bar{\lambda} = h$. We need to show that there
exists a matrix $G$ for which $\max_{h \not = 0} \frac{\| h\|_{G
Q^{-1} G^T}^2}{\| h\|_W^2} = 1$. Since we know that the maximal
value in the previous optimization problem is always smaller than
$1$ (otherwise \eqref{eq_descent} would not hold), we have to show
that there exists matrix $G$ and vector $h$ for which:
\begin{align}
\label{tight}
 \| h\|_{G Q^{-1} G^T}^2
= \|h\|_W^2.
\end{align}
 Let us consider matrices $G$ with $0/\sqrt{\sigma_i}$ entries and
that have exactly $\omega$ nonzeros on each row and on each column.
If we let $h$ be the vector with all entries equal to $1$, then
\eqref{tight} holds.

We  denote by $\Lambda^*$ the set of optimal solutions of dual
problem \eqref{eq_dual_prob}. According to \cite[Theorem
2.3.2]{HirLem:96}, if Assumption \ref{ass_strong} holds for our
original problem \eqref{eq_prob_princc}, then  $\Lambda^*$ is
nonempty, convex and bounded. For any  $\lambda \in \rset^{p+q}$,
we can define the following finite quantity:
\begin{equation}
\label{eq_multipleirs_bounded}
\mathcal{R}(\lambda)=\min\limits_{\lambda^* \in \Lambda^*}
\|\lambda^*-\lambda\|_W.
\end{equation}

In this paper  we propose a  distributed dual gradient method for
which we are interested  in deriving  estimates for both dual and
primal suboptimality and also for primal feasibility violation, i.e.
for a given accuracy $\epsilon$ find a primal-dual pair
$\left(\hat{z},\hat{\lambda}\right)$ such that:
\begin{eqnarray}
\label{condition_outer} && \|\left[G \hat{z}
-g\right]_{\mathbb{D}}\|_{W^{-1}}
 \leq \mathcal{O} (\epsilon), \;\;\;  \| \hat{z}- z^* \|^2 \leq
\mathcal{O} (\epsilon), \nonumber\\
[-1.5ex]\\[-1.5ex]
&&-\mathcal{O}(\epsilon)\!\leq \!f(\hat{z}) -\! f^* \!\!\leq\!
\mathcal{O} (\epsilon) ~\text{and} ~f^*-d(\hat{\lambda}) \leq
\mathcal{O} (\epsilon). \nonumber
\end{eqnarray}


\subsection{Motivation: Distributed MPC (DMPC) for networked systems}
\label{subsec_motivation}

We consider a discrete time networked system, modeled  by a directed
graph $\mathcal{G}=\left(V,E\right)$, for which the set
$V=\left\{1,\dots,M\right\}$ represents the subsystems and the
adjacency matrix $E$ indicates the dynamic couplings between these
subsystems. Note that in these settings, the graph $\mathcal{G}$ is
a particular case of the bipartite graph previously presented for
which we have $V_1=V_2=V$. The dynamics of the subsystems can be
defined by the following linear state equations \cite{NecNed:11}:
\begin{equation}
\label{mod3} x_i(t+1) =\sum _{j \in \mathcal {N}_{i} }
\bar{A}_{ij} x_j(t)+
 \bar{B}_{ij} u_j(t)  \qquad \forall i \in V,
\end{equation}
where $x_i(t) \in \mathbb{R}^{n_{x_i}}$ and $u_i(t) \in
\mathbb{R}^{n_{u_i}}$ represent the state and  the input of $i$th
subsystem at time $t$, $\bar{A}_{ij} \in \mathbb{R}^{n_{x_i} \times
n_{x_j}}$ and $\bar{B}_{ij} \in \mathbb{R}^{n_{x_i} \times
n_{u_j}}$. Note that in this case $\mathcal {N}_{j}$ denotes
the set of subsystems, including $j$, whose dynamics directly affect
the dynamics of subsystem $j$ and $\bar{\Ncal}_i$ represents the set
of subsystems, including $i$, whose dynamics are affected by the
dynamics of subsystem $i$. We also impose coupled state and input
constraints:
\begin{align}
\label{mpc_constraints}
 \sum_{j \in \Ncal_i } \bar{C}_{ij} x_j(t)
+\tilde{C}_{ij} u_j(t) \leq c_i \qquad \forall i \in V, \;\; t \geq
0.
\end{align}
For a prediction horizon of length $N$, we consider
strongly convex stage and final costs for each subsystem: $
\sum_{t=0}^{N-1} \ell_i(x_i(t),u_i(t)) +
\ell_i^{\mathrm{f}}(x_i(N))$, where the final costs
$\ell_i^{\mathrm{f}}$ and the terminal sets $X_i^{\mathrm{f}}$ are
chosen  such that the control scheme ensures closed-loop stability
\cite{MayRaw:00}.

The centralized MPC problem for the networked system \eqref{mod3},
for a given initial state $x=\left[x_1^{T}\cdots x_M^{T}\right]^T$,
can be posed as the following separable convex optimization problem
\cite{NecNed:11}:
\begin{align}
\label{mod5}
 &  \min _{x_i(t),u_i(t)}  \sum_{i=1}^M \sum_{t=0}^{N-1}
\ell_i(x_i(t),u_i(t)) + \ell_i^{\mathrm{f}}(x_i(N))  \\
& \text{s.t.:} \;\; x_i(t+1) =
\sum _{j \in \mathcal {N}^{i} } \bar{A}_{ij}  x_j (t) + \bar{B}_{ij} u_j(t), \; x_i(0) =x_i,  \nonumber \\
& \;\; \sum_{j \in \Ncal_i} \bar{C}_{ij} x_j(t) + \tilde{C}_{ij}
u_j(t) \leq \bar c_i, \; x_i(N) \in X_i^{\mathrm{f}}  ~  \forall i
\in V, ~ \forall t. \nonumber
\end{align}

\noindent For the state and input trajectory of subsystem $i$ and
the overall state and input trajectory we use the notations:
\begin{align*}
z_i \!=\!\! \left[u_{i}(0)^T \! x_i(1)^T \!\cdots\! u_{i}(N\!-\!1)^T
\! x_i(N)^T \right]^T\!\!\!,\ z\!\!=\!\!\left[z_{1}^T \!\cdots\!
z_{M} ^T\right]^T
\end{align*}
and for the total local cost over the prediction horizon
\begin{align*}
&f_i(z_i)=\sum_{t=0}^{N-1} \ell_i(x_i(t),u_i(t)) +
\ell_i^{\mathrm{f}}(x_i(N)).
\end{align*}
In these settings,  the centralized MPC problem for
the networked system \eqref{mod3}, for a given initial state
$x=\left[x_1^{T}\cdots x_M^{T}\right]^T$, can be posed as the
separable convex optimization problem \eqref{eq_prob_princc}, where
$n_i=N(n_{u_i}+n_{x_i})$, the equality constraints $Az=b$ are
obtained by stacking all the dynamics \eqref{mod3} together, while
the inequality constraints $Cz \leq c$ are obtained by writing the
state and input constraints \eqref{mpc_constraints} in compact form,
over the prediction horizon (see e.g. \cite{NecSuy:09}). Note also
that for the matrices $A$ and $C$, each block $A_{ji}$ and $C_{ji}$
is zero when $E_{ji}=0$.

In the following sections, we analyze the structural  properties of the dual problem
\eqref{eq_dual_prob} and then we propose a fully distributed dual
gradient method for solving this problem which exploits the
separability of the dual function and allow us to recover a suboptimal
and nearly feasible solution for our original problem
\eqref{eq_prob_princc} in linear time.


\section{Error bound property of the dual problem}
\label{sec_error_bound}

In this section, under Assumption \ref{ass_strong}, we prove an
error bound type property on the corresponding dual problem
\eqref{eq_dual_prob}. For completeness, first  we briefly review the
existing results on error bound properties for a convex problem in
the form:  \[ \min_{y \in D} \psi(y), \] where $\psi(\cdot)$ is
convex function, with Lipschitz continuous gradient, and $D$ is a
polyhedral set.  We are interested in finding optimal points for
this problem, i.e. points $y$ satisfying $y = [y - \nabla
\psi(y)]_D$.  Typically, in order to show linear convergence for
gradient based methods used for solving the above problem,  we need
to require some nondegeneracy assumption on the problem (e.g. strong
convexity) which does not hold for many practical applications (e.g.
(DMPC), (NUM) or (DC-OPF) problems). A new line of analysis, that
circumvents these difficulties, was developed using the notion of
error bound, which estimates the distance to the solution set from
an $y \in D$ by the norm of the proximal residual $\nabla^+ \psi(y)
= [y - \nabla \psi(y)]_D - y$ (in \cite{FacPan:03} $\nabla^+
\psi(y)$ is referred to as the \textit{natural map}). For objective
functions of the form $\psi(y) = \bar{\psi}(G^Ty)$, with
$\bar{\psi}(\cdot)$ strongly convex function and  $G$  a general
matrix, the authors in \cite{LuoTse:93}  show a \textit{local} error
bound property that holds in a neighborhood of the solution set,
while in \cite{WanLin:13} the authors show a global error bound
property provided that the set $D$ is a bounded polyhedron or the
entire space.

Our approach for proving  a global error bound property for the dual
problem \eqref{eq_dual_prob} is in a way similar to the one in
\cite{LuoTse:93,LuoTse:93a,WanLin:13}.  However, our results are
more general in the sense that:  we derive a \textit{global} error
bound property as opposed to the results in
\cite{LuoTse:93,LuoTse:93a} where the authors show this property
only \textit{locally} in a neighborhood of the solution set, and we
allow the constraints set to be an \textit{unbounded} polyhedron, as
opposed to the results in \cite{WanLin:13} where the authors  show
an error bound property only for constraints defined in terms of
\textit{bounded} polyhedra. Also, our proximal residual introduced
below is more general than the one used in the standard analysis of
the error bound property (see e.g.
\cite{LuoTse:93,LuoTse:93a,WanLin:13}). Last but
not least important is that our approach also works for dual
problems, which allows us to prove for the first time a global error
bound property for such problems.

For the convex function $f$, we denote its conjugate
\cite{RocWet:98}:
\begin{equation*}
\tilde{f}(y)=\sum_{i=1}^M \tilde{f}_i(y),
\end{equation*}
where $\tilde{f}_i(y) =\max \limits_{z_i \in \rset^{n_i}} \langle y,
z_i \rangle - f_i(z_i)$. According to Proposition 12.60 in
\cite{RocWet:98}, under the Assumption \ref{ass_strong} $(a)$ (in
particular, under the assumption that $f_i$ has Lipschitz gradient),
each function $\tilde{f}_i(y)$ is strongly convex w.r.t. Euclidean
norm, with constant $\frac{1}{L_i}$, which implies that function
$\tilde{f}$ is strongly convex w.r.t. the same norm, with constant:
\[\sigma_{\tilde{\mathrm{f}}}=\sum_{i=1}^M \frac{1}L_i.\] Note that in these settings
our dual function can be written as:
\begin{equation}
\label{eq_dual_conj} d(\lambda)=-\tilde{f}(-G^T \lambda)-g^T
\lambda.
\end{equation}
Note that if $G$ has full row rank (and thus $p+q \leq n$),  then it
follows immediately that the dual function $d$ is strongly concave.
Therefore, we consider below the nontrivial case when  $p+q > n$,
i.e. $G$ has not full row rank. Recall that for the projection of
$z$ onto the set $\mathbb{D}$ w.r.t. the norm $\|\cdot\|_W$ we use
$\left[z\right]_{\mathbb{D}}^W$. We denote further the proximal
residual:
\begin{equation}
\label{eq_gradient_map} \nabla^+
d(\lambda)=\left[\lambda+W^{-1}\nabla
d(\lambda)\right]_{\mathbb{D}}^W-\lambda ~~~\forall \lambda \in
\mathbb{D}.
\end{equation}

The following lemma whose proof can be found e.g. in \cite[Lemma
6.4]{NecCli:13} (see also \cite{WanLin:13, LuoTse:93})   will help
us prove the desired error bound property for our dual problem
\eqref{eq_dual_prob}.
\begin{lemma}
\label{lemma_unique_tstar}  Let Assumption \ref{ass_strong} hold.
Then,  there exists a unique $y^* \in \rset^n$ such that:
\begin{equation}
\label{eq_unique_tstar} G^T \lambda^* =y^* \quad \forall \lambda^*
\in \Lambda^*.
\end{equation}
Moreover, $\nabla d(\lambda) = G \nabla \tilde{f}(-y^*)-g$ is
constant for all $\lambda \in \Lambda$, where the set
$\Lambda=\left\{\lambda \in \mathbb{D} : G^T \lambda =
y^*\right\}$.
\end{lemma}

The following theorem, which is one of the main results of the
paper, establishes a global error bound like property for our dual
problem \eqref{eq_dual_prob}:
\begin{theorem}
\label{theorem_error_bound} Let Assumption \ref{ass_strong} hold.
Then,  there exists a constant $\kappa$, depending on  the data of
problem \eqref{eq_prob_princc} and the term
$\mathcal{T}(\lambda) = \max \limits_{\lambda^*
\in \Lambda^*}\|\lambda - \lambda^*\|_W$, such that the following
error bound property holds for dual problem
\eqref{eq_dual_prob}:
\begin{equation}
\label{eq_error_bound} \|\lambda - \bar{\lambda}\|_W \leq
\kappa\left( \mathcal{T}(\lambda) \right) \|\nabla^+
d(\lambda)\|_W ~~ \forall \lambda \in \mathbb{D},
\end{equation}
where $\bar{\lambda}=\left[\lambda\right]_{\Lambda^*}^W$ and $\kappa \left( \mathcal{T}(\lambda) \right)$ is given by $\kappa \left( \mathcal{T}(\lambda) \right)
=\theta_1^2 \frac{4}{\sigma_{\tilde{\mathrm{f}}}} +
12\theta_2^2\!\left(2 \mathcal{T}^2(\lambda) \!+\!2\|\nabla
d(\bar{\lambda})\|_W^2\!\right)\!\left(\!1\!+\!3\theta_1^2\frac{2}{\sigma_{\tilde{\mathrm{f}}}}\right)$,
with  $\theta_1$ and $\theta_2$ positive constants depending
on problem~data.
\end{theorem}
\pf Since the proof is involved and makes use of some technical
results, for clarity of the exposition we present it in the
Appendix. \qed

\noindent In general, it is difficult to derive good estimates for
the constant $\kappa$ which depends on the Hoffman's bound for
polyhedra \cite{Rob:73}. However, there are  special classes of
optimization problems when $\kappa$ can be computed explicitly:
e.g., if matrix $G$ has full row rank, then $d$ is $\sigma_{d,W}$ -
strongly concave w.r.t. the norm $\|\cdot\|_W$ and $\kappa$ was
already computed in Pang \cite{Pan:87} as
$\kappa=\frac{2}{\sigma_{d,W}}$; for other special cases see e.g.
\cite{NecCli:13,WanLin:13}.

\noindent Based on Theorem \ref{theorem_error_bound} we will prove
in the following section the linear rate of convergence of a
distributed dual gradient method. To our knowledge this  is the
first result showing  \textit{global} linear convergence rate on
primal suboptimality and infeasibility for the last primal iterate
of a dual gradient algorithm, as opposed e.g.  to the results in
\cite{LuoTse:93a} where only \textit{local} linear convergence was
derived for such an algorithm or results in
\cite{BecNed:13,GisDoa:13,KogFin:11,NecNed:13,NedOzd:09,PatBem:12}
where \textit{sublinear} convergence is proved.


\section{Linear convergence for dual distributed gradient
 method under an error bound property}
\label{sec_dg_error_bound}

The  existing convergence results from the literature on dual
gradient methods   either  show sublinear rate of convergence
\cite{BecNed:13,GisDoa:13,KogFin:11,NecNed:13,NedOzd:09,PatBem:12}
or at most \textit{local} linear convergence \cite{LuoTse:93a}. In
this section we show, that under the error bound property for the
dual problem as proved in Theorem \ref{theorem_error_bound}, which
is valid for quite general assumptions (see Assumption
\ref{ass_strong}),  we have linear convergence for a distributed
dual gradient method. Thus, we now introduce the following  fully
distributed dual gradient method:
\begin{center}
\framebox{
\parbox{7cm}{
\begin{center}
\textbf{ Algorithm {\bf (DG)}}
\end{center}
{Initialization: $\lambda^0 \in \mathbb{D}$. For $k\geq 0$ compute:}
\begin{enumerate}
\item{$z^k = \arg \min\limits_{z \in \rset^n}
\mathcal{L}(z,\lambda^k)$}.

\item ${\lambda}^{k+1} = \left[\lambda^k+W^{-1}\nabla
d(\lambda^k)\right]_{\mathbb{D}}$.
\end{enumerate}
}}
\end{center}
Note that  if we cannot solve the inner problem (step 1) exactly,
but with some inner accuracy, then our framework allows us  to use
approximate solutions $z^k$ and inexact dual gradients. This is
beyond the scope of the present paper, but for more details see e.g.
\cite{NecNed:13}. The main difference between our  Algorithm
(\textbf{DG}) and the algorithms proposed in literature
\cite{DoaKev:11,GisDoa:13,NecNed:13,NecSuy:08,NedOzd:09,MeiUlb:13}
consists in the way we update the sequence $\lambda^k$. Instead of
using a classical projected gradient step with a scalar centralized
step size as in
\cite{DoaKev:11,GisDoa:13,NecNed:13,NecSuy:08,NedOzd:09,MeiUlb:13},
we update $\lambda^k$ using a projected weighted gradient step which
allows us to obtain a fully distributed scheme. The following
relation, which is a generalization of a standard result for
gradient methods shows that Algorithm (\textbf{DG}) is an ascent
method \cite{Nes:04}:
\begin{equation}
\label{eq_ascent_dual} d(\lambda^{k+1}) \geq
d(\lambda^k)+\frac{1}{2}\|\lambda^k-\lambda^{k+1}\|_W^2  \quad  \forall
k \geq 0.
\end{equation}
Using further Lemma \ref{descent_lemma} with $\lambda=\lambda^1$ and
$\bar{\lambda}=\lambda^0$ we have:
\begin{align*}
d(\lambda^1) &\geq d(\lambda^0) + \langle \nabla d(\lambda^0), \lambda^1-\lambda^0 \rangle -\frac{1}{2}\|\lambda^1-\lambda^0\|_W^2\\
&=\max \limits_{\lambda \in \mathbb{D}} d(\lambda^0) + \langle \nabla d(\lambda^0), \lambda-\lambda^0 \rangle -\frac{1}{2}\|\lambda-\lambda^0\|_W^2\\
&\geq \max \limits_{\lambda \in \mathbb{D}} d(\lambda) -\frac{1}{2}\|\lambda-\lambda^0\|_W^2 \geq f^* -\frac{1}{2}\mathcal{R}^2(\lambda^0),
\end{align*}
where the first  equality follows from the definition of $\lambda^1$
in Algorithm (\textbf{DG}), the second inequality from the concavity
of function $d$ and the last inequality from
$\Lambda^* \subseteq \mathbb{D}$ and definition of $\max$. Using
now the previous relation we obtain:
\begin{equation}
\label{eq_upper_lambda1}
f^*-d(\lambda^1) \leq \frac{1}{2}\mathcal{R}^2(\lambda^0).
\end{equation}
The next lemma will help us
analyze the convergence of the Algorithm (\textbf{DG}):

\begin{lemma}
\label{lemma_decrese_dist_dg} Let Assumption \ref{ass_strong}
hold and the sequence $\{\lambda^k\}_{k \geq 0}$ be generated by
Algorithm (\textbf{DG}). Then, the following inequalities hold:
\begin{equation}
\label{eq_decrease} \!\|\lambda^k-\lambda^*\|_W \!\leq\! \cdots \leq
\|\lambda^0-\lambda^*\|_W \quad \forall \lambda^* \in \Lambda^*, k
\geq 0.
\end{equation}
\end{lemma}
\pf First we notice that the update
${\lambda}^{k+1}$ can be also viewed as the unique optimal solution of the maximization
of the following quadratic approximation of $d$:
\begin{equation}
\label{eq_projection_dual} \max \limits_{\lambda
\in \mathbb{D}} d(\lambda^k)+\langle \nabla d(\lambda^k),
\lambda-\lambda^k \rangle -\frac{1}{2}\|\lambda-\lambda^k\|_W^2.
\end{equation}
Taking now $\lambda=\lambda^*$ in the optimality condition of
\eqref{eq_projection_dual}, we obtain the following inequality:
\begin{equation}
\label{eq_opt_cond_decrease} \langle \nabla d(\lambda^k)-W
(\lambda^{k+1}-\lambda^k),\lambda^*-\lambda^{k+1}\rangle \leq 0.
\end{equation}
Further, we can write:
\begin{align}
\label{inequalities_gradient}
&\|\lambda^{k+1}\!\!-\!\lambda^*\|_W^2\!
=\!\|\lambda^{k+1}-\lambda^k+\lambda^k-\lambda^*\|_W^2 \\
&=\!\|\lambda^k\!\!-\!\lambda^*\|_W^2\!+\!2\langle W(
\lambda^{k+1}\!\!-\!\lambda^k),\lambda^{k}\!\!-\!\lambda^{k+1}\!\!+\!\lambda^{k+1}\!\!-\!\lambda^*
\rangle\! \nonumber\\
&\qquad \qquad+\!\|\lambda^{k+1}\!\!-\!\lambda^k\|_W^2 \nonumber\\
&=\!\|\lambda^k\!\!-\!\!\lambda^*\|^2_W\!+\!2\langle W(
\lambda^{k+1}\!\!-\!\!\lambda^k),\!\lambda^{k+1}\!\!-\!\!\lambda^*
\rangle\!-\!\|\lambda^{k+1}\!\!-\!\!\lambda^k\|_W^2 \nonumber \\
&\!\leq \|\lambda^k-\lambda^*\|^2_W-2\langle \nabla
d(\lambda^k),\lambda^*-\lambda^k\rangle \nonumber \\
&\qquad+2\left(\langle \nabla
d(\lambda^k),\lambda^{k+1}-\lambda^k\rangle\!-\frac{1}{2}
\|\lambda^{k+1}\!\!-\lambda^k\|^2_W\right)  \nonumber\\
&\!\leq
\|\lambda^k-\lambda^*\|^2_W+\!2\left(d(\lambda^k)\!-\!d(\lambda^*)\right)\!+\!2\left(d(\lambda^{k+1})-d(\lambda^k)\right)
\nonumber \\
&\!=\|\lambda^k-\lambda^*\|^2_W+2\left(d(\lambda^{k+1})-d(\lambda^*)\right)
\leq \|\lambda^k-\lambda^*\|^2_W,\nonumber
\end{align}
where the first inequality follows from \eqref{eq_opt_cond_decrease}
and the second one is derived from the concavity of the function $d$
and Lemma \ref{descent_lemma}. \qed

Using now inequality \eqref{eq_decrease} in \eqref{eq_error_bound}
we obtain one important relation that estimates the
distance from the dual optimal solution set of the sequence $\lambda^k$ in terms of
the norm of a proximal residual:
\begin{equation}
\label{global_error_bound} \|\lambda^k-\bar{\lambda}^k\|_W \leq
\overline{\kappa} \;  \|\nabla^+ d(\lambda^k)\|_W \quad \forall k
\geq 0,
\end{equation}
where from Theorem \ref{theorem_error_bound} we have:
\begin{align*}
\overline{\kappa}= \theta_1^2
\frac{4}{\sigma_{\tilde{\mathrm{f}}}}\!+\!12\theta_2^2\left(2\mathcal{T}^2(\lambda^0)\!+\!2
\mathcal{T}^2_1(\Lambda^*)\right)\left(1\!+\!3\theta_1^2
\frac{2}{\sigma_{\tilde{\mathrm{f}}}}\right),
\end{align*}
where we have defined the positive constant $\mathcal{T}(\lambda^0)
= \max\limits_{\lambda^* \in \Lambda^*} \| \lambda^0 -
\lambda^*\|_W$, which is finite since $\Lambda^*$ is a
bounded set. Moreover,  from Lemma  \ref{lemma_unique_tstar} we have that $\nabla
d(\bar{\lambda})$ is constant for  all $\bar{\lambda} \in \Lambda^*$ and thus
we can define the positive constant $\mathcal{T}_1(\Lambda^*) =\|\nabla
d(\bar{\lambda})\|_W$  for all $\bar{\lambda} \in \Lambda^*$. Further,
since $W$ is a positive definite diagonal matrix, the following
relation is straightforward:
\begin{equation}
\label{eq_decr_map} \|\nabla^+ d(\lambda^k)\|_W
=\|\lambda^{k+1}-\lambda^k\|_W.
\end{equation}
Combining now \eqref{global_error_bound} with \eqref{eq_decr_map},
we can write:
\begin{equation}
\label{eq_relation1} \|\lambda^k-\bar{\lambda}^k\|_W \leq
\overline{\kappa}\|\nabla^+
d(\lambda^k)\|_W=\overline{\kappa}\|\lambda^{k+1}-\lambda^k\|_W.
\end{equation}
The following theorem provides an estimate on the dual
suboptimality for Algorithm ({\bf DG}) and follows similar lines as in \cite{WanLin:13,LuoTse:93a,NecCli:13}:

\begin{theorem}
\label{theorem_dual_optim_dg} Let Assumption \ref{ass_strong}  hold
and  sequences $\left(z^k,\lambda^k\right)_{k\geq 0}$ be generated
by Algorithm ({\bf DG}). Then, an estimate on dual suboptimality for
\eqref{eq_dual_prob} is given by:
\begin{equation}
\label{bound_dual_optim_dg} f^*-d({\lambda}^{k+1})\leq \frac{1}{2}
\left(\frac{4(1+\overline{\kappa})}{1+4(1+\overline{\kappa})}\right)^{k-1}
\mathcal{R}^2(\lambda^0).
\end{equation}
\end{theorem}

\pf From the optimality conditions of problem
\eqref{eq_projection_dual} we have:
\begin{equation}
\label{eq_optimality_lambdak} \langle \nabla
d(\lambda^k),\bar{\lambda}^k\!-\!\lambda^{k+1}\rangle \!\leq\! \langle W
(\lambda^{k+1}\!-\!\lambda^k), \bar{\lambda}^k\!-\!\lambda^{k+1}\rangle
\!\leq \! 0,
\end{equation}
where we recall that $\bar{\lambda}^k=[\lambda^k]^W_{\Lambda^*}$.
Further, since the optimal value of the dual function is unique we
can write:
\begin{align}
\label{ineq_linear_dg}
&f^*\!-\!d(\lambda^{k+1})=d(\bar{\lambda}^k)\!-\!d(\lambda^{k+1})\leq
\langle \nabla d(\lambda^{k+1}), \bar{\lambda}^k\!-\!\lambda^{k+1}
\rangle \nonumber\\
& =\langle \nabla d(\lambda^{k+1})\!-\!\nabla d(\lambda^k),
\bar{\lambda}^k\!-\!\lambda^{k+1}\rangle \!+\!\langle \nabla
d(\lambda^{k}), \bar{\lambda}^k\!\!-\!\lambda^{k+1} \rangle \nonumber\\
& \leq \| \nabla d(\lambda^{k+1})-\nabla d(\lambda^k)\|_{W^{-1}}
\|\bar{\lambda}^k - \lambda^{k+1}\|_{W} \nonumber\\
&\qquad \qquad \qquad+ \langle W (\lambda^{k+1}-\lambda^k),
\bar{\lambda}^k-\lambda^{k+1}\rangle \nonumber\\
& \leq \|\lambda^{k+1}-\lambda^k\|_{W} \|\bar{\lambda}^k -
\lambda^{k+1}\|_{W} \nonumber\\
&\qquad \qquad+ \|\lambda^{k+1}-\lambda^k\|_{W}
\|\bar{\lambda}^k-\lambda^{k+1}\|_{W} \nonumber\\
&=2\|\lambda^{k+1}-\lambda^k\|_{W}\|\bar{\lambda}^k-\lambda^{k+1}\|_{W},
\end{align}
where the second inequality follows from
\eqref{eq_optimality_lambdak}. Using now relation
\eqref{eq_relation1} we can write:
\begin{align*}
\|\bar{\lambda}^k-\lambda^{k+1}\|_{W} &\leq
\|\bar{\lambda}^k-\lambda^{k}\|_{W}+\|{\lambda}^k-\lambda^{k+1}\|_{W}\\
&\leq \left(1+\overline{\kappa}\right)\|{\lambda}^k-\lambda^{k+1}\|_{W}.
\end{align*}
Introducing now the previous inequality in \eqref{ineq_linear_dg}
and using \eqref{eq_ascent_dual} we have:
\begin{align*}
f^*-d(\lambda^{k+1})&\leq
2\left(1+\overline{\kappa}\right)\|{\lambda}^k-\lambda^{k+1}\|_{W}^2\\
&\leq 4\left(1+\overline{\kappa}\right)
\left(d(\lambda^{k+1})-d(\lambda^k)\right).
\end{align*}
Rearranging the terms in the previous inequality we obtain:
\begin{equation}
\label{eq_recursion} f^*-d(\lambda^{k+1}) \leq
\frac{4(1+\overline{\kappa})}{1+4(1+\overline{\kappa})}\left(f^*-d(\lambda^k)\right).
\end{equation}
Applying now \eqref{eq_recursion} recursively and using \eqref{eq_upper_lambda1} we obtain
\eqref{bound_dual_optim_dg}. \qed

The following theorems give estimates on the primal feasibility
violation and suboptimality for Algorithm (\textbf{DG}). Note that
usually, for recovering an approximate primal solution from dual
gradient based methods,  we need to use averaging (see e.g.
\cite{DoaKev:11,KogFin:11,NecNed:13,NecSuy:08,NedOzd:09,PatBem:12}).
In what follows, we do not consider averaging and we prove  linear
convergence for  the last primal iterate, a result which appears to
be new.
\begin{theorem}
\label{theorem_primal_fesa_dg} Under the assumptions of Theorem
\ref{theorem_dual_optim_dg}, the following estimate holds for the
primal infeasibility:
\begin{equation*}
\left\|\left[G z^k\!-\!g
\right]_{\mathbb{D}}\right\|_{W^{-1}} \!\leq\!
\left(\frac{4(1+\overline{\kappa})}{1+4(1+\overline{\kappa})}\right)^{\frac{k-2}{2}}\mathcal{R}(\lambda^0).
\end{equation*}
\end{theorem}
\pf Using the descent property of dual gradient method
\eqref{eq_ascent_dual} we have:
\begin{align}
\label{eq_feasibility1_dg} \|\lambda^{k}\!-\!\lambda^{k+1}\|_W^2 &\leq
2 \left(d(\lambda^{k+1})-d(\lambda^k)\right) \leq
2\left(f^*\!-\!d(\lambda^k)\right) \nonumber\\
&\leq\!
\!\left(\!\frac{4(1+\overline{\kappa})}{1\!+\!4(1\!+\!\overline{\kappa})}\!\right)^{k-2}\mathcal{R}^2(\lambda^0),
\end{align}
where in the last inequality we used Theorem
\ref{theorem_dual_optim_dg}. In order to prove the statement of the theorem we will first show that $\left\|\left[\nabla
d(\lambda^{k})\right]_{\mathbb{D}}\right\|_{W^{-1}}^2 \leq
\|\lambda^{k}\!\!-\!\lambda^{k+1}\|_W^2$. We will prove this
inequality componentwise. First, we recall that
$\mathbb{D}=\rset^p \times \rset^q_+$. Thus, for all $i=1,\dots,p$
we have:
\begin{align}
\label{eq_part1} &\left|\left[\nabla_i
d(\lambda^{k})\right]_{\rset}\right|_{W_{ii}^{-1}}^2\!\!=\left|\nabla_i
d(\lambda^{k})\right|_{W_{ii}^{-1}}^2 \\
&~~=\left|\lambda^{k}_i-\lambda^{k}_i-W_{ii}^{-1}\nabla_i
d(\lambda^{k})\right|_{W_{ii}}^2 =\left|\lambda^{k}_i-\lambda^{k+1}_i\right|_{W_{ii}}^2,
\nonumber
\end{align}
where in the last equality we used the definition of
$\lambda^{k+1}$. We now introduce  the following disjoint sets:
$I_{-}=\left\{i \in [p+1,p+q] : \nabla_i d(\lambda^{k}) < 0
\right\}$ and $I_{+}=\left\{i \in [p+1,p+q] : \nabla_i
d(\lambda^{k}) \geq 0 \right\}$. Using these notations and the
definition of $\mathbb{D}$, we can write:
\begin{equation}
\label{eq_part2} \left|\left[\nabla_i
d(\lambda^{k})\right]_{\rset_+}\right|_{W_{ii}^{-1}}^2 \!\!=\! 0
\leq \left| \lambda^{k}_i-\lambda^{k+1}_i\right|_{W_{ii}}^2  ~
\forall i \in I_{-}.
\end{equation}
On the other hand, for all $i \in I_{+}$ we have:
\begin{align}
\label{eq_part3} &\left|\!\left[\nabla_i
d(\lambda^{k})\right]_{\rset_+}\right|_{W_{ii}^{-1}}^2\!\!=\!\left|\nabla_i
d(\lambda^{k})\right|_{W_{ii}^{-1}}^2 \\
&\qquad \quad=\left|\!\left[W_{ii}^{-1}\nabla_i
d(\lambda^{k})\right]_{\rset_+}\right|_{W_{ii}}^2 =\left|\lambda^{k}_i-\lambda^{k+1}_i\right|_{W_{ii}}^2.
\nonumber
\end{align}
Summing up the relations \eqref{eq_part1},\eqref{eq_part2} and
\eqref{eq_part3} for all $i=1,\dots,p+q$ we obtain:
\begin{equation}
\label{eq_grad_equiv}
\left\|\left[\nabla
d(\lambda^{k})\right]_{\mathbb{D}}\right\|_{W^{-1}}^2 \leq
\left\|\lambda^{k}-\lambda^{k+1}\right\|_{W}^2.
\end{equation}
Combining now \eqref{eq_grad_equiv} and \eqref{eq_feasibility1_dg} we obtain:
\begin{align*}
\left\|\left[\nabla
d(\lambda^{k})\right]_{\mathbb{D}}\right\|_{W^{-1}}^2 &\leq
\|\lambda^{k}-\lambda^{k+1}\|_W^2\\
&\leq \left(\frac{4(1+\overline{\kappa})}{1+4(1+\overline{\kappa})}\right)^{k-2}\mathcal{R}^2(\lambda^0).
\end{align*}
Taking into account the definition of $\nabla  d$ we conclude the
statement. \qed

We now characterize  the primal suboptimality and the distance from
the last iterate $z^{k}$, generated by Algorithm (\textbf{DG}), to
the optimal solution $z^*$.
\begin{theorem}
\label{theorem_primal_optim_dg} Under the assumptions of Theorem
\ref{theorem_primal_fesa_dg}, the following estimate on primal
suboptimality for problem \eqref{eq_prob_princc} can be derived:
\begin{align}
\label{bound_primal_optim_dg}
-\left(\frac{4(1+\overline{\kappa})}{1+4(1+\overline{\kappa})}\right)^{\frac{k-2}{2}}
&\left(\mathcal{R}^2(\lambda^0)+\mathcal{R}(\lambda^0)\|\lambda^0\|_{W}\right)
\nonumber \\
& \leq f(z^k)-f^* \leq v(k),
\end{align}
where
\begin{align*}
v(k)&=\frac{\left(\mathcal{R}^2(\lambda^0)+ \mathcal{R}(\lambda^0)\|\lambda^0\|_W\right)}{\underline{w}\sqrt{\sigma_{\text{f}}}}\|G\|\!\!\left(\frac{4(1+\overline{\kappa})}{1+4(1+\overline{\kappa})}\right)^{\frac{k-2}{2}}\\
&\qquad+\frac{\max\limits_{i=1,\dots,M}L_i}{2\sigma_{\text{f}}}
\left(\frac{4(1+\overline{\kappa})}{1+4(1+\overline{\kappa})}\right)^{k-2}\mathcal{R}^2(\lambda^0),
\end{align*}
with $\underline{w}=\lambda_{\text{min}}(W)$. Moreover, the sequence
$z^{k}$ converges to the unique optimal solution $z^*$ of
\eqref{eq_prob_princc} with  the following rate:
\begin{equation}
\label{bound_primal_solution_dg} \|z^{k}-z^*\| \leq \sqrt{\frac{1}{\sigma_{\text{f}}}}
\left(\frac{4(1+\overline{\kappa})}{1+4(1+\overline{\kappa})}\right)^{\frac{k-2}{2}}
\mathcal{R}(\lambda^0).
\end{equation}
\end{theorem}
\pf In order to prove the left-hand side inequality of
\eqref{bound_primal_optim_dg} we can write:
\begin{align}
\label{ineq_aux_optim_dg} f^*& =d(\lambda^*)=\min_{z \in \rset^n}
f(z)+\langle
\lambda^*,Gz-g\rangle \\
&\leq f(z^k)+\langle \lambda^*,G z^k-g \rangle \nonumber \\
&\leq f(z^k)+
\langle \lambda^*,\left[G z^k-g\right]_{\mathbb{D}} \rangle \nonumber \\
&\leq f(z^k) + \|\lambda^*\|_{W}\left\|\left[G
z^k-g\right]_{\mathbb{D}}\right\|_{W^{-1}} \nonumber\\
& \leq f(z^k)+
\left(\mathcal{R}(\lambda^0)+\|\lambda^0\|_{W}\right)\left\|\left[G
z^k-g\right]_{\mathbb{D}}\right\|_{W^{-1}}, \nonumber
\end{align}
where the second inequality follows from the fact that $\lambda^*
\in \mathbb{D}$ and the third one from the Cauchy-Schwartz inequality.
Using now Theorem \ref{theorem_primal_fesa_dg} we obtain the result.

For proving the right hand-side
inequality of \eqref{bound_primal_optim_dg} we first show
\eqref{bound_primal_solution_dg}. First, let us note that since
 $f$ is $\sigma_{\mathrm{f}}$-strongly convex w.r.t. Euclidean norm, it
follows that $\lu(z,\lambda^k)$ is also
$\sigma_{\mathrm{f}}$-strongly convex in the variable $z$ in the
same norm. We recall that $d(\lambda^k)=f(z^k) + \langle \lambda^k,
G z^k - g \rangle$ and $\nabla d(\lambda^k)=Gz^k-g$. Taking now into
account that $z^k =\arg \min_{z \in \rset^n} \lu(z,\lambda^k)$ and
the fact that $\langle \lambda^k, \nabla d(\lambda^*) \rangle \leq
0$, from the strong convexity of $\lu$ we have:
\begin{align*}
\frac{\sigma_{\mathrm{f}}}{2}\|z^k\!-\!z^*\|^2\!&\leq\!
\lu(z^*,\lambda^k)\!-\!\lu(z^k,\lambda^k)\\
&= \! f(z^*)\!+\!\langle \lambda^k,\!\!\nabla
d(\lambda^*)\rangle\!-\!\!f(z^k)\!-\!\langle \lambda^k,\!\nabla
d(\lambda^k)\rangle\\
&= f^* \!+\!\langle \lambda^k,\!\!\nabla
d(\lambda^*)\rangle\! - d(\lambda^k)  \\
&\leq f^*-d(\lambda^k).
\end{align*}
Using further Theorem \ref{theorem_dual_optim_dg} in the previous inequality we obtain \eqref{bound_primal_solution_dg}. From  the
Lipschitz continuity property of $\nabla f$ we obtain:
\begin{align*}
f(z^{k})-f^* &\leq \langle \nabla f(z^*), z^k-z^*
\rangle+\frac{\max_i L_i}{2}\|z^{k}-z^*\|^2\\
&=\langle -G^T \lambda^*, z^k-z^* \rangle+\frac{\max_i L_i}{2}\|z^{k}-z^*\|^2\\
&\leq\! \frac{\|\lambda^*\|_W \|G\|}{\underline{w}}\|z^{k} \!-\!
z^*\| \!+\! \frac{\max_i L_i}{2}\|z^{k} \!-\! z^*\|^2,
\end{align*}
where the first equality is deduced from the optimality conditions
of problem $z^*=\arg\min f(z)+\langle \lambda^*, Gz-g\rangle$
and in the last inequality we used the
Cauchy-Schwartz inequality and the fact that $\|\cdot\| \leq
\frac{1}{\underline{w}}\|\cdot\|_W$ and $\|G^T \lambda^*\| \leq
\|G\| \|\lambda^*\|$. Using now the definition of
$\mathcal{R}(\lambda^0)$ and
\eqref{bound_primal_solution_dg} we obtain the result. \qed \\



\section{Distributed implementation}
\label{sec_distributed}

In this section we analyze the distributed implementation of
Algorithm (\textbf{DG}). We look first at step $1$ of the algorithm.
According to \eqref{eq_separable_dual}, for all $i \in V_1$ we have:
\begin{align}
\label{distributed_primal_update} z_i^{k}&=\arg\min_{z_i \in
\rset^{n_i}} f_i(z_i) +\left\langle \lambda^k, \left[A_{ i}^T
C_{ i}^T\right]^Tz_i \right\rangle \nonumber\\
&=\arg \min_{z_i \in \rset^{n_i}} f_i(z_i) +\sum_{j \in \bar{\Ncal}_i}
\left(\left[A_{ji}^T C_{ji}^T\right]\lambda_j^k\right)^T z_i.
\end{align}
Thus, in order to compute $z_i^k$ the algorithm requires only
local information, namely
$\left\{A_{ji},C_{ji},\lambda^k_j\right\}_{j \in \bar{\Ncal}_i}$. Using now the definitions of $W$ and
$\nabla d$, step $2$ in Algorithm (\textbf{DG}) can be written in
the following form for all $j \in V_2$:
\begin{equation}
\label{distributed_dual_update}
\lambda_j^{k+1}=\left[\lambda_j^k+\left[\begin{array}{c}
   W_{\nu,jj}^{-1} \sum_{i \in \Ncal_j} A_{ji} z_i^k\\
   W_{\mu,jj}^{-1} \sum_{i \in \Ncal_j} C_{ji} z_i^k\\
\end{array}\right]\right]_{\rset^{p_j}\times \rset^{q_j}_+},
\end{equation}
where $W_{\nu,jj}$ and $W_{\mu,jj}$ denote the $j$th
block-diagonal element of matrix $W_{\nu}$ and $W_{\mu}$,
respectively. Taking into account the definitions of $W_{\nu,jj}$
and $W_{\mu,jj}$ we can conclude that in order to update the dual
variable $\lambda_j^{k+1}$ in step $2$ of Algorithm
(\textbf{DG}) we require only local information
$\left\{L_{\text{d}_i},A_{ji},C_{ji},z_i^k\right\}_{i \in
\Ncal_j}$.

We discuss further some implementation issues for the case of (DMPC) problems.
A standard approach for such problems is to consider quadratic stage and final costs for each subsystem, i.e.:
\begin{align*}
\ell_i(x_i(t),u_i(t))&=0.5 \|x_i(t)\|^2_{\bar Q_i}+0.5\|u_i(t)\|^2_{\bar R_i} ~\text{and}\\
\ell_i^{\text{f}}(x_i(N))&=0.5\|x_i(N)\|^2_{\bar P_i},
\end{align*}
where $\bar Q_i$, $\bar R_i$ and $\bar P_i$ are positive definite matrices of appropriate dimensions. In this case,
each objective function is quadratic, i.e. $f_i(z_i)=0.5 \|z_i\|^2_{Q_i}$, where $Q_i$ is
given by: $Q_i=diag\left(I_{N-1}\otimes\left[
                                    \begin{array}{cc}
                                      {\bar R}_i & 0 \\
                                      0 & {\bar Q}_i \\
                                    \end{array}
                                  \right], {\bar R}_i, {\bar P}_i
\right)$.  Further, note that the matrices $A$ and $C$ are block-sparse having each block $(i,j)$ zero for every subsystem $i$ which is not influenced
by subsystem $j$. According to \eqref{distributed_primal_update}, the update $z_i^k$ of subsystem $i$ can be done in closed form:
\begin{equation*}
z_i^{k}= Q_i^{-1} \sum_{j \in \bar{\Ncal}_i}
\left(\left[A_{ji}^T C_{ji}^T\right]\lambda_j^k\right),
\end{equation*}
and requires only information from those subsystems that are
influenced by subsystem $i$, i.e. $\bar{\Ncal}_i$. Similarly, for
updating $\lambda_j^{k+1}$ corresponding to subsystem $j$, we need
information from its neighbors, i.e. $\Ncal_j$  (those subsystems
that affect subsystem $j$). Moreover,   if we define  the
\textit{measure of sparsity} for  the incidence matrix $E$ as: \[
\omega=\max_{i,j \in V} \left\{ |\bar{\Ncal}_i|, |\Ncal_j|\right\},
\] then due to the structure of block matrices $Q_i, A_{ji}$ and
$C_{ji}$, the computational complexity of one step of Algorithm
(\textbf{DG}) for (DMPC) is linear in both, the number of subsystems
$M$ and the horizon length $N$, i.e. $\mathcal{O}(MN)$, provided
that $\omega \ll M$. Finally, our algorithm is scalable in the sense
that removing or adding a new node (subsystem) can be done
immediately using only local information.

Further, we note that all the estimates for the convergence rate
for primal and dual suboptimality and primal feasibility violation
derived in Section \ref{sec_dg_error_bound} depends on the upper
bound on the norm of the optimal Lagrange multipliers
$\mathcal{R}$, which at its turn depends on the degree of
separability of problem \eqref{eq_prob_princc}, characterized by
the sets $\Ncal_i$ and $\bar{\Ncal}_j$. In order to see this
dependence we can write further:
\begin{equation*}
\mathcal{R}^2(\lambda^0)\!=\!\!\max_{\lambda^* \in
\Lambda^*}\!\|\lambda^*\!-\lambda^0\|_W^2=\!\!\max_{\lambda^* \in \Lambda^*}
\sum_{j=1}^{\bar{M}} \sum_{i \in \bar{\Ncal}_j}\!L_{\mathrm{d}_i}
\|\lambda^*_j-\lambda_j^0\|^2,
\end{equation*}
from which it is straightforward to notice that $\mathcal{R}$
depends on the cardinality of each $\bar{\Ncal}_j$. On the other
hand, for each $i$ we recall that:
\begin{equation*}
L_{\text{d}_i}=\frac{\left\|\left[
\begin{array}{c}
  \left[A_{ji}\right]_{j \in
\Ncal_i} \\
   \left[C_{ji}\right]_{j \in
\Ncal_i} \\
\end{array}
\right]\right\|^2}{\sigma_i},
\end{equation*}
which depends on the cardinality of the set $\Ncal_i$. Thus, we
can conclude that $\mathcal{R}$ depends on the cardinality of
$\Ncal_i$ and $\bar{\Ncal}_j$ which represent a natural measure
for the degree of separability of our original problem
\eqref{eq_prob_princc}.


\section{Numerical simulations}
\label{sec_numerical}

\begin{figure*}[]
\subfloat[Primal and dual suboptimality along iterations. \label{fig:iterations_a}]{%
\includegraphics[width=0.45\textwidth]{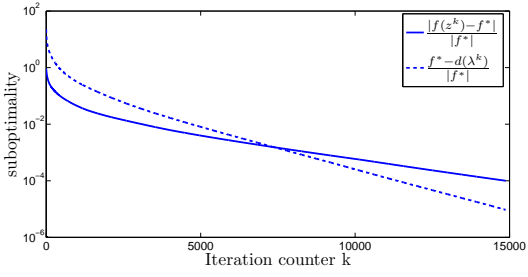}
}
\hfill
\subfloat[Ratio $k_{DG}/k_{CG}$ w.r.t. sparsity measure $\omega$. \label{fig:iterations_b}]{%
\includegraphics[width=0.45\textwidth]{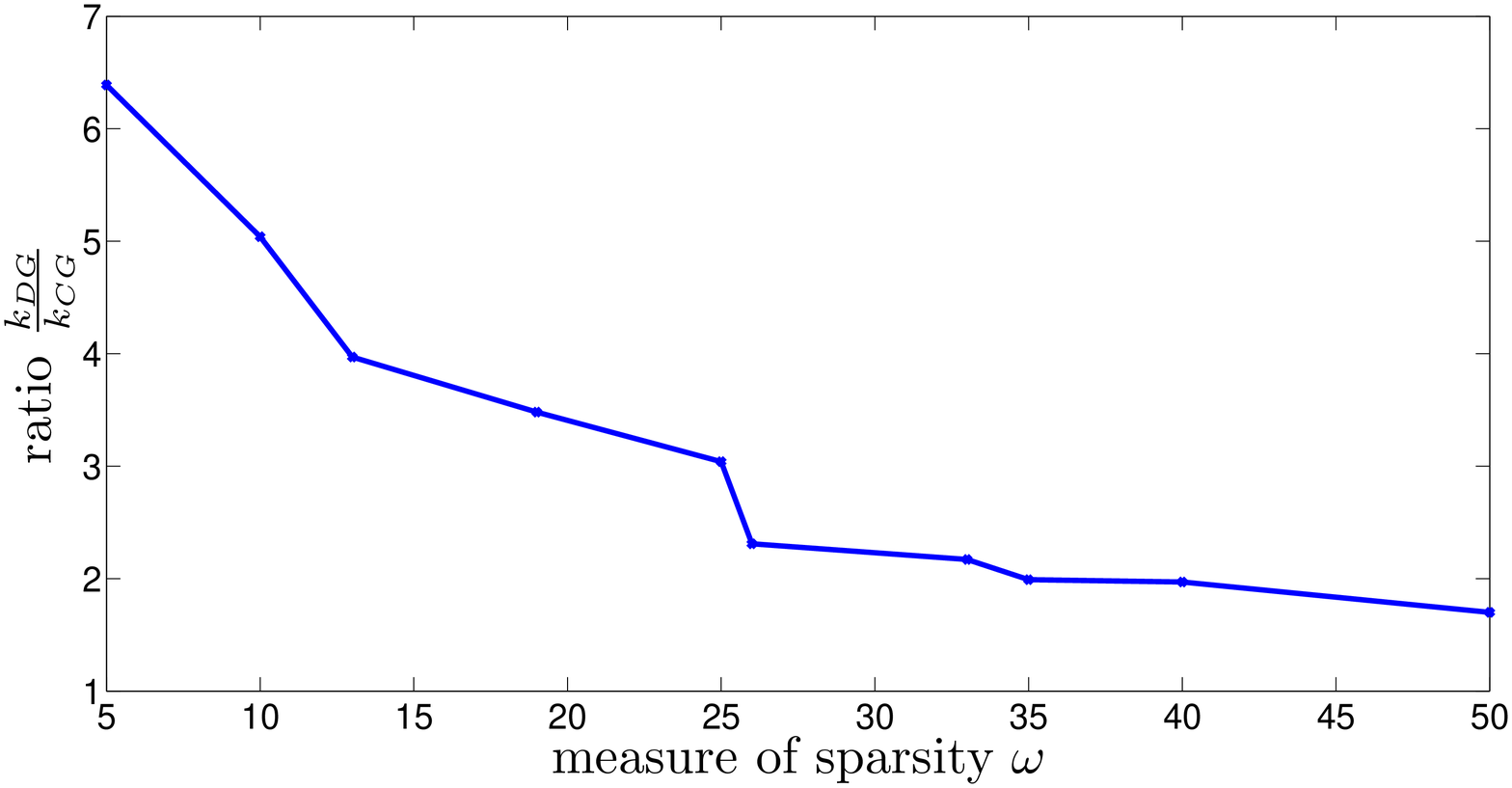}
} \caption{Behavior of Algorithm (\textbf{DG}).}
\label{fig:iterations}
\end{figure*}

In this section we consider problems of the form
\eqref{eq_prob_princc}, where the objective functions are given by:
\[ f_i(z_i) = 0.5 \|z_i\|_{Q_i}^2 + q_i^T z_i + \gamma_i \log
\left(1+e^{\langle a_i, z_i\rangle}\right). \] Note that this type
of function satisfies Assumption \ref{ass_strong} (a), provided
that $Q_i$ are positive definite matrices, and are intensively used
in (DMPC) applications for $\gamma_i=0$ or (NUM) applications for
$\gamma_i>0$. We first generate a sparse communication graph
$\mathcal{G}=\left(V,E\right)$ characterized by an incidence matrix
$E \in \rset^{M \times M}$ generated randomly with different degrees
of sparsity given by $\bar{\Ncal}_i$ and $\Ncal_j$. Recall that we
have defined  the measure of sparsity of the incidence matrix $E$
as: $\omega=\max_{i,j \in V} \left\{ |\bar{\Ncal}_i|,
|\Ncal_j|\right\}$. We take $n_i=n_j$ for all $i, j$. Matrices $Q_i
\in \rset^{n_i \times n_i}$, $A_{ji} \in
\rset^{\lceil\frac{3n_i}{4}\rceil \times n_i}$ and $C_{ji} \in
\rset^{\lceil\frac{3n_i}{2}\rceil \times n_i}$ are taken from a
normal distribution with zero mean and unit variance. Matrix $Q_i$
is then made positive definite by transformations $Q_i \leftarrow
Q_i^T Q_i+\sigma_i I_{n_i}$, where $\sigma_i$ are randomly generated
from the interval $[1,10]$. Further, vectors $b$, $c$ are chosen
such that the problem is feasible, and $q_i$ and $a_i$ are taken
from an uniform distribution and $\gamma_i=1$ for all $i$.

\begin{figure*}[]
\centerline{ \subfloat[Linear convergence of Algorithm \textbf{(DG)}
in the last iterate: logarithmic scale of primal suboptimality and
infeasibility. We also compare with the theoretical sublinear
estimates (dot lines) for the convergence rate of order ${\mathcal
O}(1/\sqrt{k})$. The plot clearly shows our theoretical findings,
i.e. linear convergence.]{
\includegraphics[width=1.0\textwidth,height=4.5cm]{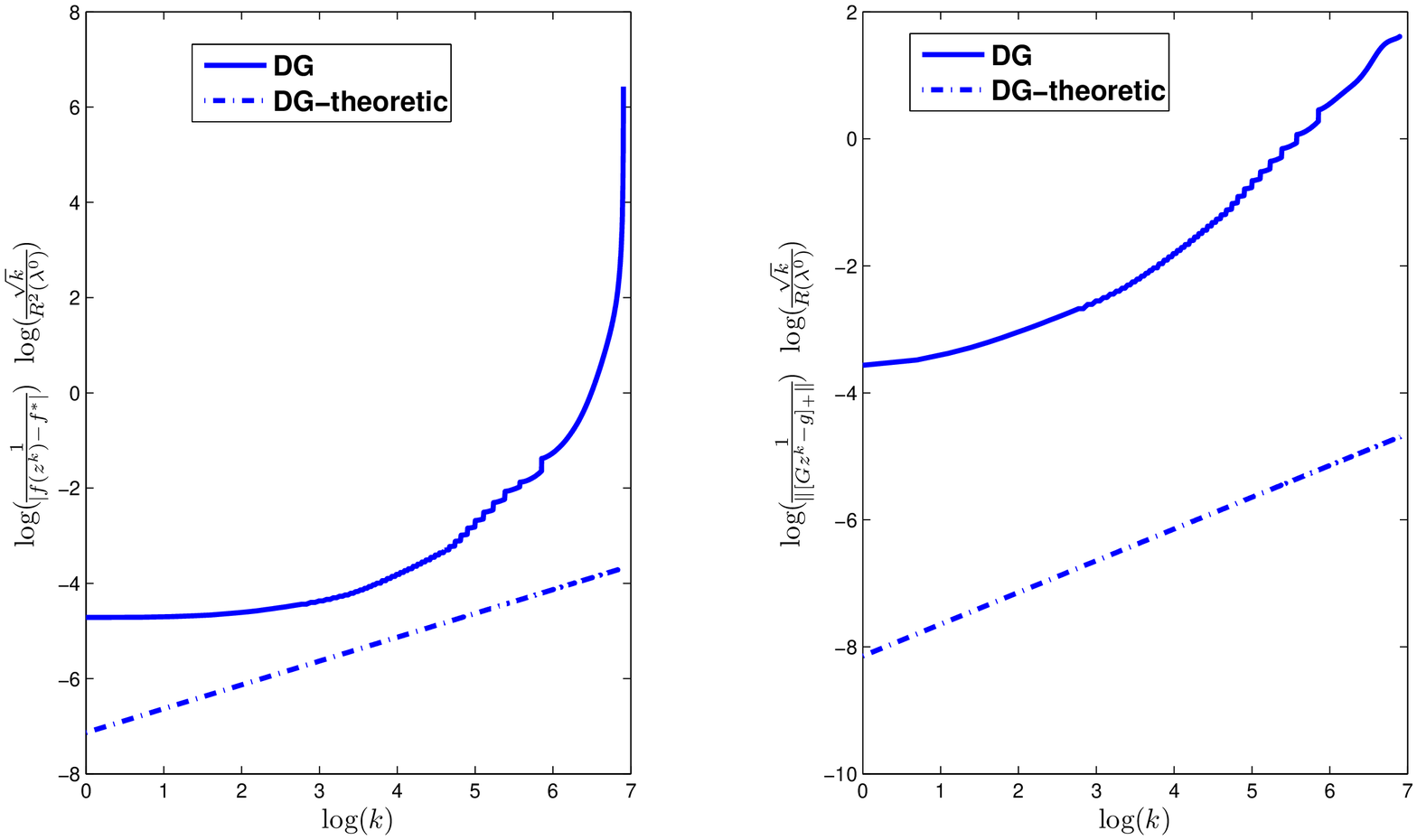}
}} \caption{Algorithm (\textbf{DG}): theoretical vs. practical
behavior.} \label{fig:t_p}
\end{figure*}

In order to analyze the behavior of Algorithm (\textbf{DG}) we first
consider a problem with the number of subsystems $M=100$, the
dimension of local variables $n_i=10$ for all $i \in
\left\{1,\dots,100\right\}$ and $\omega=15$. We are interested in
analyzing the evolution of both, primal and dual suboptimality,
w.r.t. the number of iterations. We consider an accuracy $\epsilon
=10^{-4}$ and impose a stopping criterion of the form:
\begin{equation}
\label{eq_numer_stopping}
\frac{|f(z^k)-f^*|}{|f^*|} \leq \epsilon,
\end{equation}
where $f^*$ is computed using CVX. We plot the results in
logarithmic scale in Figure \ref{fig:iterations_a}. We can observe
that both dual and primal suboptimality converge linearly, which
confirms the theoretical results derived in Section
\ref{sec_dg_error_bound}. Moreover, in one of our recent papers
\cite{NecPat:14} we have proved that dual gradient algorithm in the
last primal iterate is converging sublinearly (order
$\mathcal{O}(1/\sqrt{k})$) in terms of primal suboptimality and
infeasibility, provided that the primal objective function $f$ is
only  strongly convex (i.e. $f$ has not Lipschitz gradient and thus
no error bound property holds for the dual problem). In particular,
we have proved that primal suboptimality is of order
$\mathcal{O}(\frac{\mathcal{R}^2(\lambda^0)}{\sqrt{k}})$ and primal
infeasibility is of order
$\mathcal{O}(\frac{\mathcal{R}(\lambda^0)}{\sqrt{k}})$. From  Figure
\ref{fig:t_p} we again observe  linear convergence. In the same
figure we also  plot the theoretical sublinear estimates for the
convergence rate of order $\mathcal{O}(1/\sqrt{k})$ for Algorithm
\textbf{(DG)} in the last iterate as described above (see
\cite{NecPat:14} for more details). The plot clearly confirms our
theoretical findings, i.e. linear convergence of Algorithm
\textbf{(DG)} in the last iterate.

We are  also interested in comparing the performances of Algorithm
(\textbf{DG}) with the ones of the centralized dual gradient method,
called Algorithm (\textbf{CG}), where for updating the dual variable
$\lambda$ we use the centralized step size
$L_{\text{d}}^{-1}I_{p+q}$ (see also \cite{NecNed:13}):
\[  (\textbf{CG}): \quad  {\lambda}^{k+1} = \left[\lambda^k + L_{\text{d}}^{-1} \nabla
d(\lambda^k)\right]_{\mathbb{D}}. \] For this purpose, we consider
a set of 10 problems with fixed dimension, $M=100$ and $n=10$,
accuracy $\epsilon=10^{-2}$ and the measure of sparsity  $\omega$
ranging from $5$ to $50$. We plot in Figure \ref{fig:iterations_b}
the ratio between the number of iterations performed by Algorithm
(\textbf{DG}) ($k_{DG}$) and by Algorithm (\textbf{CG}) ($k_{CG}$),
respectively. On the one hand, we observe that for small values of
the sparsity measure $\omega$, Algorithm (\textbf{DG}) clearly
outperforms Algorithm (\textbf{CG}). On the other hand, by
increasing the sparsity measure $\omega$ we observe  a reduction in
the ratio between the number of iterations performed by the two
algorithms.

\vspace*{-0.3cm}

Further, we consider   problems of different dimensions, having the
measure of sparsity $\omega=0.15 M$, and compare the Algorithms
(\textbf{DG}) and (\textbf{CG}) in terms of number of iterations
performed for obtaining a suboptimal solution with accuracy
$\epsilon=10^{-2}$. In Table \ref{table_iterations} we present the
average results obtained for 10 randomly generated problems for each
dimension. We use the notation
$\overline{w}=\lambda_{\text{max}}(W)$ and we recall that
$\underline{w}=\lambda_{\text{min}}(W)$. We can observe that for all
dimensions, Algorithm (\textbf{DG}) outperforms  Algorithm
(\textbf{CG}) due to the reduced level of sparsity $\omega$ and
smaller values of the entries of weighted step size $W$ compared to
the centralized one $L_{\text{d}}$.
\begin{table}[ht!]
\centering {  
\begin{tabular}{|c|c|c|c|}
\hline
\backslashbox{results}{$(M,n_i)$} & $(200,10)$ & $(100,20)$ & (50,40) \\
\hline
 $k_{DG}$ &  3861  &  4117  &  4936  \\
\hline
 $k_{CG}$ &  16121  &  19541  &  27973  \\
\hline
 $L_{\text{d}}$ &  691  & 802   & 1346   \\
\hline
 $\overline{w}$ &  653  &  854  &  1433  \\
\hline
 $\underline{w}$ &  21  & 49  &  79  \\
\hline
\end{tabular}}
\begin{center}{\caption{Average number of iterations  for finding an $\epsilon$-solution.} \label{table_iterations}}\end{center}
\end{table}

\section{Conclusions}
\label{sec_conclusions}

In this paper we have proposed and analyzed   a fully distributed
dual gradient method for solving the Lagrangian dual of a primal
separable convex optimization problem with linear constraints. Under
the strong convexity and Lipschitz continuous gradient property
assumptions of the primal objective function we have provided a
global error bound for the dual problem. Using this property, we
have proved global linear convergence rate for both primal and dual
suboptimality and for primal feasibility violation for our
distributed dual gradient algorithm. We have also discussed
distributed implementation aspects for our method and provided
several numerical simulations which confirm the theoretical results
and the efficiency of our approach.


\section*{Appendix}
In order to prove Theorem \ref{theorem_error_bound} we first
need some technical results. First we recall from Lemma \ref{lemma_unique_tstar} that
there exists a unique $y^* \in
\rset^n$ such that:
\begin{equation*}
G^T \lambda^* =y^* \quad \forall \lambda^*
\in \Lambda^*.
\end{equation*}
Moreover, $\nabla d(\lambda) = G \nabla \tilde{f}(-y^*)-g$ is
constant for all $\lambda \in \Lambda$, where we have defined the set:
$$ \Lambda = \left\{\lambda \in \mathbb{D} : G^T \lambda =
y^*\right\}. $$
We introduce further the following notations:
\begin{equation}
\label{eq_notations_lambda} r=\left[\lambda\right]_{\Lambda}^W ~\text{and}
~\bar{r}=\left[r\right]_{\Lambda^*}^W \qquad \forall \lambda \in \mathbb{D}.
\end{equation}
Using now the notations \eqref{eq_notations_lambda} and $\bar{\lambda}=\left[ \lambda \right]_{\Lambda^*}^W$ we can write:
\vspace*{-0.2cm}
\begin{align}
\label{eq_intermediate_ineq}  \|\lambda- \bar{\lambda}\|^2_W & \!\leq
\left\|\lambda-\bar{r}\right\|^2_W \leq \left(\|\lambda
-r\|_W+\|r-\bar{r}\|_W\right)^2 \nonumber\\
& \!\leq 2\|\lambda-r\|_W^2 \!+ 2\|r-\bar{r}\|^2_W \quad \forall \lambda \in \mathbb{D}.
\end{align}
In what follows we show how we can find upper bounds on $\|\lambda
-r\|_W$ and $\|r-\bar{r}\|_W$ such that we will be able to establish
the error bound property on the dual problem \eqref{eq_dual_prob}
given in Theorem \ref{theorem_error_bound}.
\begin{lemma}
\label{lemma_gradient-gradient_map} Let Assumption \ref{ass_strong}
hold and $\nabla^+ d$  given in \eqref{eq_gradient_map}. Then, the
following inequality holds for all $\lambda, \omega \in \mathbb{D}$:
\begin{equation*}
\langle \nabla d(\omega)-\nabla d(\lambda),\!\lambda-\omega \rangle\!\!
\leq \!2 \|\nabla^+ \!d(\lambda)-\nabla^+ \!d(\omega)\|_W \;
\|\lambda-\omega\|_W.
\end{equation*}
\end{lemma}

\pf First, let us recall that $\left[\lambda+W^{-1}\nabla
d(\lambda)\right]^W_{\mathbb{D}}$ is the unique solution of the
optimization problem:
\begin{equation}
\label{eq_projection_W} \min_{\xi \in \mathbb{D}}
\|\xi-\lambda-W^{-1}\nabla d(\lambda)\|_W^2,
\end{equation}
for which the optimality conditions reads:
\begin{align*}
&\!\!\!\!\!\!\!\left\langle W\!\left(\!
\left[\lambda\!+\!W^{-1}\nabla
d(\lambda)\right]^W_{\mathbb{D}}\!\!-\!\left(\!\lambda\!+\!W^{-1}\nabla
d(\lambda)\right)\!\!\right),\right. \\
&\qquad \qquad \left.\xi -\left[\lambda\!+\!W^{-1}\nabla
d(\lambda)\right]^W_{\mathbb{D}} \right\rangle \geq 0 \quad \forall \xi
\in \mathbb{D}.
\end{align*}
Taking now $\xi=\left[\omega+W^{-1}\nabla d(\omega
)\right]^W_{\mathbb{D}}$ in the previous  inequality, adding and
subtracting both, $\lambda$ and $\omega$, in the right term of the
scalar product and using the definition of $\nabla^+ d$ we obtain:
\begin{equation*}
\langle W\!\left(\nabla^+\! d(\lambda)\!-\!W^{-1}\nabla d(\lambda)\right)\!,\!
\nabla^+\! d(\lambda)+\lambda - \omega -\nabla^+\! d(\omega) \rangle \!\!\leq \!0,
\end{equation*}
and the symmetry of matrix $W$ leads to:
\begin{align*}
&\langle \nabla^+ d(\lambda)-W^{-1}\nabla
d(\lambda),\\
&\qquad \qquad ~~W\left(\lambda-\omega\right)+ W\left(\nabla^+
d(\lambda)-\nabla^+ d(\omega)\right) \rangle  \leq 0.
\end{align*}
Rearranging  the terms in the previous inequality we get:
\begin{align*}
-&\langle \nabla d(\lambda),\lambda-\omega\rangle\\
&\leq \!\!-\langle
\nabla^+ d(\lambda), W\!(\lambda\!-\!\omega) \rangle \!+\! \langle \nabla
d(\lambda), \nabla^+ d(\lambda)\!-\!\nabla^+ d(\omega)\rangle\\
&\qquad \qquad - \langle \nabla^+ d(\lambda), W \left(\nabla^+
d(\lambda)- \nabla^+ d(\omega)\right) \rangle.
\end{align*}
Writing now the previous inequality with $\lambda$ and $\xi$
interchanged and summing them up we can write:
\begin{align*}
&\langle \nabla d(\xi)\!-\!\nabla d(\lambda), \lambda\!-\!\xi
\rangle\\
&\leq \!\langle
\nabla^+\!d(\xi)\!-\!\nabla^+\!d(\lambda), W(\lambda\!-\!\xi)
\rangle\!\\
&\qquad+\langle \nabla d(\lambda)\!-\!\nabla
d(\xi),
\nabla^+\!d(\lambda)\!-\!\nabla^+\!d(\xi)\rangle\!\\
&\qquad-\!\|\nabla^+d(\lambda)\!-\!\nabla^+
d(\xi)\|_W^2\\
&\leq \langle \nabla^+\!d(\xi)\!-\!\nabla^+\!d(\lambda),
W(\lambda\!-\!\xi) \rangle\!\\
&\qquad+\langle \nabla d(\lambda)\!-\!\nabla
d(\xi), \nabla^+\!d(\lambda)-\nabla^+ d(\xi) \rangle\\
&\leq \|\nabla^+ d(\lambda)-\nabla^+ d(\xi)
\|_W\|W(\lambda-\xi)\|_{W^{-1}}\\
&\qquad +\|\nabla^+ d(\lambda)-\nabla^+ d(\xi)
\|_W\|\nabla d(\lambda)-\nabla
d(\xi)\|_{W^{-1}}\\
&\leq 2\|\nabla^+ d(\lambda)-\nabla^+ d(\xi) \|_W
\|\lambda-\xi\|_W,
\end{align*}
which concludes the statement. \qed

The next lemma gives un upper bound on $\|\lambda -r\|_W$:
\begin{lemma}
\label{lemma_upper_r} Under Assumption \ref{ass_strong} there exists
a constant $\kappa_1$ such that the following inequality holds:
\begin{equation}
\label{eq_bound_r} \|\lambda-r\|^2_W \leq
 \kappa_1 \; \|\nabla^+d(\lambda)\|_W  \; \|\lambda-\bar{\lambda}\|_W \quad \forall \lambda
\in \mathbb{D},
\end{equation}
where $\kappa_1=\frac{2}{\sigma_{\tilde{\mathrm{f}}}}\theta_1^2$ and
 $\theta_1$  depends on the matrix $G$.
\end{lemma}
\pf  First, let us notice that we can write the set $\Lambda$
explicitly as:
\begin{align}
\label{eq_sets_mangasarian2} \Lambda&=\left\{\omega \in
\rset^{p+q}:~F \omega \leq 0,~~ G^T \omega = y^* \right\},
\end{align}
where $F=\left[0_{q,p} ~ -I_q\right]$. Since $\lambda \in
\mathbb{D}$, it implies $F \lambda \leq 0$ and therefore, according
to Theorem 2 in \cite{Rob:73}, we can bound the distance between a
vector $\lambda$ and the polyhedron $\Lambda$ as follows:
\begin{equation}
\label{eq_mangasarian} \|\lambda-r\|_W \leq \theta_1 \|G^T \lambda -
y^*\|_{\infty} \leq \theta_1 \|G^T \lambda - y^*\|,
\end{equation}
where $\theta_1$ is the Hoffman's bound depending  on the matrix $G$
and on the norms $\|\cdot\|_W$ and $\|\cdot\|_{\infty}$ (see eq. (6)
in \cite{Rob:73} for a formula  to compute Hoffman's bound).  From
the strong convexity property of $\tilde{f}$ combined with the fact
that $G^T \bar{\lambda}=y^*$ we have:
\begin{align}
\label{eq_strong_fstar_bound} &\sigma_{\tilde{\mathrm{f}}}\|G^T
\lambda -y^*\|^2 \\
&\leq \langle \nabla \tilde{f}(-G^T
\lambda)-\nabla \tilde{f}(-G^T
\bar{\lambda}), -G^T \lambda +G^T \bar{\lambda} \rangle \nonumber\\
&=\langle -G\nabla \tilde{f}(-G^T \lambda)+g +G\nabla
\tilde{f}(-G^T
\bar{\lambda})-g, \lambda - \bar{\lambda} \rangle  \nonumber \\
&=\langle \nabla d(\bar{\lambda}) -\nabla d(\lambda),
\lambda-\bar{\lambda}\rangle \nonumber \\
&\leq 2\|\nabla^+ d(\lambda) - \nabla^+
d(\bar{\lambda})\|_W\|\lambda-\bar{\lambda}\|_W \nonumber\\
&=2\|\nabla^+ d(\lambda)\|_W\|\lambda-\bar{\lambda}\|_W, \nonumber
\end{align}
where the last inequality follows from Lemma
\ref{lemma_gradient-gradient_map}  and the last equality follows  from the fact
that $\nabla^+d(\bar{\lambda})=0$ for $\bar{\lambda} \in \Lambda^*$.
Combining now \eqref{eq_mangasarian} with
\eqref{eq_strong_fstar_bound} we obtain the result. \qed

The next result establishes an upper bound on $\|r-\bar{r}\|_W$:
\begin{lemma}
\label{lemma_upper_ru} Let Assumption \ref{ass_strong} be satisfied. Then, the following inequality is valid:
\begin{equation}
\label{eq_upper_ru} \|r -  \bar{r}\|_W^2 \!\leq\! \kappa_2(\mathcal{T}(\lambda)) \|\nabla^+\!
d(\lambda)\|_W\|\lambda - \bar{\lambda}\|_W \!\!\! \quad \forall \lambda \!\in\!
\mathbb{D},
\end{equation}
where $r, \bar{r}$ are given in \eqref{eq_notations_lambda},  $\bar{\lambda}=\left[ \lambda \right]_{\Lambda^*}^W$ and
\[\kappa_2(\mathcal{T}(\lambda)) \!= 6\theta_2^2 \left( 2 \mathcal{T}^2(\lambda) +2\|\nabla
d(\bar{\lambda})\|_{W^{-1}}^2 \right) \! \left(1 \!+\!
3\theta_1^2\frac{2}{\sigma_{\tilde{\mathrm{f}}}}\right),\] with
$\mathcal{T}(\lambda) = \max \limits_{\lambda^* \in
\Lambda^*}\|\lambda-\lambda^*\|_W$ and   $\theta_2$ being a constant
depending on $C$, $\nabla d(\bar{\lambda})$ and $y^*$.
\end{lemma}
\pf Since $\Lambda^* \subseteq \Lambda \subseteq \mathbb{D}$ and
$G^T \xi = y^*$ for all $\xi \in \Lambda$,  the dual problem
\eqref{eq_dual_prob} has the same optimal solutions as the following
linear program:
\begin{align}
\label{eq_linear1} \arg \max_{\lambda \in \mathbb{D}}
d(\lambda)&\!=\!\arg \max_{\xi \in \Lambda} d(\xi)\!=\! \arg \max_{\xi \in
\Lambda} -\tilde{f}(-y^*) -\langle g,\xi\rangle \nonumber \\
&\!=\!\arg \max_{\xi \in \Lambda} -\langle g,\xi\rangle.
\end{align}
Further, let us recall that $\nabla d(\zeta)= G \nabla
\tilde{f}(-y^*)-g$ for any $\zeta \in \Lambda$ and thus we
have that $\langle \nabla d(\zeta), \xi\rangle= \langle
\nabla \tilde{f}(-y^*), y^*\rangle-\langle g, \xi \rangle$ for all
$\zeta, \xi \in \Lambda$. Therefore, we can write further for any $\zeta \in \Lambda$:
\begin{align}
\label{eq_linear2} \arg \max_{\xi \in \Lambda} \langle \nabla
d(\zeta), \xi\rangle &=\arg \max_{\xi \in \Lambda}  \langle
\nabla \tilde{f}(-y^*), y^*\rangle -\langle g, \xi\rangle \nonumber \\
&=\arg \max_{\xi \in \Lambda} -\langle g, \xi\rangle.
\end{align}
Combining now \eqref{eq_linear1} with \eqref{eq_linear2}, we can
conclude that any solution of the dual problem \eqref{eq_dual_prob}
$\bar{\xi}=[\xi]_{\Lambda^*}^W$, with $\xi \in \Lambda$,  is also a
solution of linear program  \eqref{eq_linear2}. Since  $\Lambda^* \subseteq \Lambda$, then
 for any $\left[ \lambda \right]_{\Lambda^*}^W = \bar{\lambda} \in \Lambda^*$ we have that   $\nabla_{\nu}
d(\bar{\lambda})=Az^*-b = 0$ and $\nabla_{\mu}
d(\bar{\lambda})=Cz^*-c \leq 0$, and thus we also have that the
maximum in \eqref{eq_linear2} is finite for any   $\zeta =
\bar{\lambda} = \left[ \lambda \right]_{\Lambda^*}^W \in \Lambda^*$.
Thus problem \eqref{eq_linear2} is solvable for  any $\zeta =
\bar{\lambda} = \left[ \lambda \right]_{\Lambda^*}^W \in \Lambda^*$.
Applying now Theorem 2 in \cite{Rob:73} to the optimality conditions
of problem \eqref{eq_linear2} and its dual we obtain:
\begin{equation}
\label{eq_robinson} \|\xi-\bar{\xi}\|_W \leq \theta_2 |\langle
\nabla d(\bar{\lambda}), \xi \rangle -\langle \nabla
d(\bar{\lambda}), \bar{\xi}\rangle| \quad \forall \xi \in \Lambda,
\end{equation}
where $\theta_2$ is Hoffman's bound depending only on the matrix $C$
and vectors $\nabla d(\bar{\lambda})$ and $y^*$ (see eq. (6) in
\cite{Rob:73} for details). Using the previous relation we have:
\begin{align}
\label{eq_linear3} \|\xi-\bar{\xi}\|_W & \leq \theta_2 |\langle
\nabla d(\bar{\lambda}), \xi \rangle -\langle \nabla
d(\bar{\lambda}), \bar{\xi}\rangle| \nonumber \\
& =\theta_2 \langle \nabla d(\bar{\lambda}), \bar{\xi}-\xi \rangle.
\end{align}
For any $\xi \in \Lambda$ the optimality conditions of the
following  projection problem $\min \limits_{\omega \in
\Lambda}\|\omega - \xi - W^{-1}\nabla d(\bar{\lambda}) \|_W^2$
become:
\begin{align*}
&\!\!\!\!\!\!\!\!\!\left\langle W\left(\left[\xi+W^{-1}\nabla
d(\bar{\lambda})\right]_{\Lambda}^W-\xi-W^{-1}\nabla
d(\bar{\lambda})\right),\right.\\
&\qquad \qquad\left.\left[\xi+W^{-1}\nabla
d(\bar{\lambda})\right]_{\Lambda}^W- \omega \right\rangle \leq 0,~~\forall \omega \in \Lambda.
\end{align*}
Taking now $\omega=
\bar{\xi}=[\xi]^W_{\Lambda^*}$ and since $W$ is a symmetric matrix
we obtain:
\begin{align*}
&\langle\nabla
d(\bar{\lambda}),\bar{\xi}\!-\!\xi\rangle\\
&\leq \left\langle
\left[\xi\!+\!W^{-1}\nabla
d(\bar{\lambda})\right]_{\Lambda}^W\!\!\!-\!\xi, \right.\\
& \quad \qquad \qquad \qquad \left. W\!\left(
\bar{\xi}\!-\!\left[\xi\!+\!W^{-1}\nabla
d(\bar{\lambda})\right]_{\Lambda}^W\right)\!+\!\nabla d(\bar{\lambda}) \right\rangle\\
&= \!\left\langle \left[\xi\!+\!W^{-1}\nabla
d(\bar{\lambda})\right]_{\Lambda}^W\!\!-\xi, \right. \\
& ~\quad \quad \left. W\!\left(
\xi\!-\!\left[\xi\!+\!W^{-1}\nabla
d(\bar{\lambda})\right]_{\Lambda}^W\right)+W(\bar{\xi}\!-\!\xi)\!+\!\nabla d(\bar{\lambda})\right\rangle\\
&\leq \langle \left[\xi+W^{-1}\nabla
d(\bar{\lambda})\right]_{\Lambda}^W-\xi,W(\bar{\xi}-\xi)+\nabla
d(\bar{\lambda})\rangle\\
&\leq \left\|\left[\xi+W^{-1}\nabla
d(\bar{\lambda})\right]_{\Lambda}^W-\xi\right\|_W\|W(\bar{\xi}-\xi)\|_{W^{-1}} \\
&\qquad + \left\|\left[\xi+W^{-1}\nabla
d(\bar{\lambda})\right]_{\Lambda}^W-\xi\right\|_W\|\nabla
d(\bar{\lambda})\|_{W^{-1}}\\
&=\|\nabla^{+} d(\xi)\|_W\left(\|\xi-\bar{\xi}\|_{W}+\|\nabla
d(\bar{\lambda})\|_{W^{-1}}\right),
\end{align*}
where in the last equality we used the definition of $\nabla^+ d$
and the fact that $\nabla d(\bar{\lambda})= \nabla d(\xi)$ for all
$\xi \in \Lambda$ (see Lemma \ref{lemma_unique_tstar}). Combining
now the previous inequality with \eqref{eq_linear3} and  taking
$\xi=r \in \Lambda$ we obtain:
\begin{equation}
\label{eq_bound_r_inter} \|r\!-\!\bar{r}\|_W^2\!\!\leq\!\!\theta_2^2\!\!
\left(\|r\!-\!\bar{r}\|_W\!\!+\!\!\|\nabla
d(\bar{\lambda})\|_{W^{-1}}\!\right)^2\!\|\nabla^+\!d(r)\|_W^2.
\end{equation}
Since $\bar{r}=\left[r\right]_{\Lambda^*}^W$ and $\Lambda^*
\subseteq \Lambda$ we also have
$\bar{r}=\left[\bar{r}\right]_{\Lambda}^W$. Thus, using the
nonexpansive  property of the projection we can bound the term $\|r-\bar{r}\|_W$ from above with a finite positive constant $\mathcal{T}(\lambda) = \max \limits_{\lambda^* \in
\Lambda^*}\|\lambda-\lambda^*\|_W$:
\begin{equation}
\label{eq_bound_r_max} \|r-\bar{r}\|_W \leq \|\lambda- \bar{r}\|_W
\leq \max \limits_{\lambda^* \in
\Lambda^*}\|\lambda-\lambda^*\|_W = \mathcal{T}(\lambda).
\end{equation}
Note that $\mathcal{T}(\lambda)$ is finite for any $\lambda \in \mathbb{D}$, provided that $\Lambda^*$ is a bounded set. Further, our goal is to find an upper bound for $\|\nabla^+
d(r)\|_W$ in terms of $\|\lambda-\bar{\lambda}\|_W$ and
$\|\nabla^+ d(\lambda)\|_W$. To this purpose, let us first prove
that $\nabla^+ d$ is Lipschitz continuous with constant $3$ w.r.t. the norm $\|\cdot\|_W$. For any $\lambda, \tilde{\lambda} \in
\mathbb{D}$ we can write:
\begin{align}
\label{eq_mapping_lipsch} &\|\nabla^+ d(\lambda)\!-\!\nabla^+
d(\tilde{\lambda})\|_W \leq \|\lambda\!-\!
\tilde{\lambda}\|_W \\
&\qquad \qquad+\left\|\left[\lambda\!+\!W^{-1}\nabla
d(\lambda)\right]_{\mathbb{D}}^W-\left[\tilde{\lambda}\!+\!W^{-1}\nabla
d(\tilde{\lambda})\right]_{\mathbb{D}}^W\right\| \nonumber\\
&\leq\|\lambda\!-\!\tilde{\lambda}\|_W\!+\!\|\lambda+W^{-1}\nabla
d(\lambda)\!-\!\tilde{\lambda}\!-\!W^{-1}\nabla d(\tilde{\lambda})\|_W \nonumber\\
&\leq 2\|\lambda - \tilde{\lambda}\|_W+\|\nabla d(\lambda)-\nabla
d(\tilde{\lambda})\|_{W^{-1}} \leq 3\|\lambda-\tilde{\lambda}\|_W.
\nonumber
\end{align}
Using now \eqref{eq_mapping_lipsch} with
$\tilde{\lambda}=\bar{\lambda}$ and taking into account that
$\nabla^+ d(\bar{\lambda})=0$, we have:
\begin{equation*}
\|\nabla^+ \!d(\lambda)\|_W\!=\!\|\nabla^+
d(\lambda)\!-\!\nabla^+ d(\bar{\lambda})\|_W \!\leq\!
3\|\lambda\!-\!\bar{\lambda}\|_W.
\end{equation*}
Using now again \eqref{eq_mapping_lipsch} and the  previous
inequality we get:
\begin{align}
\label{eq_ineq_error_bound}  \|\nabla^+ \!d(r)\|_W^2 &\!\leq
\left(\|\nabla^+ d(\lambda)\|_W\!+\!\|\nabla^+
d(r)\!-\!\nabla^+ d(\lambda)\|_W\right)^2 \nonumber \\
& \!\leq 2\|\nabla^+ d(\lambda)\|_W^2+2\|\nabla^+ d(r)-\nabla^+
d(\lambda)\|_W^2 \nonumber\\
&\!\leq 6\|\nabla^+
d(\lambda)\|_W\|\lambda-\bar{\lambda}\|_W+18\|\lambda-r\|_W^2 \nonumber\\
&\!\leq
6\!\left[1\!+\!3\theta_1^2\frac{2}{\sigma_{\tilde{\mathrm{f}}}}\right]\!\|\nabla^+\!
d(\lambda)\|_W\|\lambda\!-\!\bar{\lambda}\|_W,
\end{align}
where in the last inequality we used \eqref{eq_bound_r}. Introducing
now \eqref{eq_bound_r_max} and \eqref{eq_ineq_error_bound} in
\eqref{eq_bound_r_inter} and using the inequality $(\alpha+\beta)^2
\leq 2\alpha^2+2\beta^2$ we obtain the result. \qed


\pf {\it (Proof of Theorem \ref{theorem_error_bound})}  The result
given in Theorem \ref{theorem_error_bound} follows immediately by
using \eqref{eq_bound_r} from Lemma \ref{lemma_upper_r} and
\eqref{eq_upper_ru} from Lemma \ref{lemma_upper_ru} in
\eqref{eq_intermediate_ineq} and dividing both sides by $\|\lambda
-\bar{\lambda}\|_W$. \qed


\bibliographystyle{plain}
\bibliography{bibliografie2013}

\end{document}